\documentclass[smallcondensed]{svjour3}
\smartqed

\usepackage{latexsym}
\usepackage{amsmath}
\usepackage{amssymb}
\usepackage{amsfonts}
\usepackage{verbatim}
\usepackage{graphicx}
\usepackage{epstopdf}
\usepackage{subfigure}
\usepackage{float}
\usepackage{multirow}
\usepackage{color}
\usepackage[misc]{ifsym}
\usepackage{color}
\usepackage{enumerate}
\usepackage{tikz}
\usepackage{caption}
\usepackage{algorithm}
\usepackage{algorithmic}
\usepackage{setspace}
\usepackage{pifont}
\usepackage{mathrsfs}

\newcommand{\circled}[1]{\lower.7ex\hbox{\tikz\draw (0pt, 0pt)%
		circle (.5em) node {\makebox[1em][c]{\small #1}};}}
\usepackage[backref,colorlinks,linkcolor=red,anchorcolor=green,citecolor=blue]{hyperref}
\renewcommand{\arraystretch}{1.1}

\let \ssection=\section
\renewcommand{\section}{\setcounter{equation}{0}\ssection}
\DeclareMathOperator*{\diag}{diag}

\def\x{{\bf x}}
\def\y{{\bf y}}

\def\m{{\bf m}}

\newcommand{\la}{\lambda}
\def\d{{\rm{d}}}

\date{\today}

\begin{document}


\title{Efficient Spectral
Methods for Quasi-Equilibrium Closure
Approximations of Symmetric Problems on Unit Circle and Sphere}

\author{Shan Jiang \and Haijun Yu}
 \institute{S.Jiang\and H. Yu (\Letter)\\
 	\email{jshan@lsec.cc.ac.cn} (S. Jiang),
 	\email{hyu@lsec.cc.ac.cn} (H. Yu)\\
     School of Mathematical Sciences,
     University of Chinese Academy of Sciences, Beijing
     100049, China.\\
 	NCMIS \& LSEC, Institute of
 	Computational Mathematics and Scientific/Engineering
 	Computing, Academy of Mathematics and Systems
 	Science, Beijing 100190, China.\\
    ORCiD: 0000-0002-5742-0327 (H. Yu)
 	}
\titlerunning{Spectral Bingham Moment Closure on Unit Circle and Sphere}
\authorrunning{S. Jiang \and H. Yu}

\maketitle

\begin{abstract}

   Quasi-equilibrium approximation is a widely used closure approximation
   approach for model reduction with applications in complex fluids, materials
   science, etc. It is based on the maximum entropy principle and leads to
   thermodynamically consistent coarse-grain models. However, its high
   computational cost is a known barrier for fast and accurate
   applications.
   Despite its
   good mathematical properties, there are very few works on the fast and
   efficient implementations of quasi-equilibrium approximations. In this paper,
   we give efficient implementations of quasi-equilibrium
   approximations for antipodally symmetric problems on unit circle and unit
   sphere using polynomial and piecewise polynomial approximations. Comparing
   to the existing methods using linear or cubic interpolations, our approach
   achieves high accuracy (double precision) with much less storage cost. The
   methods proposed in this paper can be directly extended to handle
   other moment closure approximation problems.


\end{abstract}


\section{Introduction}

Model reduction is a classical method to obtain
computable low-dimensional mathematical models for complex systems.
Famous examples include the reduction from the Schr\"odinger
equation\cite{schrodinger_undulatory_1926} to density
function theory\cite{kohn_selfconsistent_1965}, Grad's thirteen moment model
for the Boltzmann equation\cite{grad_kinetic_1949a}, etc.
Model reduction also plays a key role in the development of polymeric
material science \cite{doi_theory_1986}, where
physically sound dynamical models are built, which are high dimensional
Fokker–Planck equations describe the evolution of $d+n$ dimensional
configuration distribution functions (CDF) of polymeric molecules. Here $d\le
3$ is the number of spatial dimensions, $n$ is the number of molecular
configurational dimensions. It is
usually impossible to solve the full $d+n$ dimensional equations
for complex systems. A common computable approach is to derive the evolution
eequations for some lower-order moments of the high-dimensional CDF from
the Fokker-Planck equation. However, except for some simple
cases (e.g. Hookean spring molecules), one usually obtain non-closed equations,
since the equations for low-order moments may involve higher-order moments. To
close these equations, one must express these higher-order
moments in terms of the lower-order moments, which is known as moment closure
problem.

The moment closure problem has been under investigation for
many years. Let's take the dynamics of liquid crystal polymer (LCP) as an
example, whose high-dimensional Fokker-Planck equation is the Doi-Smoluchowski
model\cite{doi_theory_1986,yu_nonhomogeneous_2010a}.
Various closure approximations for this model have been proposed, such as the
Doi's
quadratic closure \cite{doi_theory_1986}, the Hinch-Leal closure
\cite{hinch_constitutive_1976},
orthotropic closure \cite{cintrajr_orthotropic_1995} and the Bingham closure
\cite{chaubal_closure_1998}.
Feng et al. \cite{feng_closure_1998} examined the performances of five
commonly
used closures by numerical
simulations and found that the Bingham closure gives best results.
In fact, the Bingham closure is a particular
case of quasi-equilibrium approximation (QEA) for antipodally symmetric CDF on
$n$-sphere
(e.g. CDF for rod-like polymers), which is
an application of the maximum entropy principle (MEP) to dynamical systems that
is widely used
in statistical physics. The earliest application of the MEP can date back to
Gibbs's classic work \cite{gibbs_elementary_1902}.
Modern applications of MEP starts from Jaynes\cite{jaynes_information_1957}.
A systemic depiction of QEA for model reduction is given by Gorban et al.
\cite{gorban_corrections_2001,gorban_constructive_2004,gorban_invariant_2005}.
In the polymeric dynamics field, Chaubal and Leal first applied
quasi-equilibrium closure approximation
to rod-like polymer systems by using Bingham
distribution\cite{chaubal_closure_1998}, which
is the maximum entropy antipodally symmetric distribution on unit sphere given
second order moments\cite{bingham_antipodally_1974}.
Ilg et al. gave a system
analysis of QEA with applications to flexible polymers in homogeneous
system\cite{ilg_canonical_2002} and rod-like polymers\cite{ilg_canonical_2003},
and proved validity of energy dissipation for homogeneous systems.
Yu et al. \cite{yu_nonhomogeneous_2010a} applied the Bingham closure to
nonhomogeneous LCP systems and developed a relatively
simple but general nonhomogeneous kinetic model
for LCPs as well as efficient reduced moment models that maintain energy dissipation.

Despite its good mathematical and physical properties, efficient numerical
implementation of QEA is not an easy task. Chaubal and Leal used a global cubic
polynomial approximation fitted by the method of least square to implement
Bingham closure\cite{chaubal_closure_1998}, which has relatively large
numerical error. Grosso et al. \cite{grosso_closure_2000} gave an efficient
implementation
of Bingham closure by using Cayley-Hamilton theorem with symmetric properties
of the moments,
where a global quadratic approximation is used, which also results in large
numerical error. Yu et al. \cite{yu_nonhomogeneous_2010a} designed an efficient
implementation of Bingham closure for $2$-dimensional problem, where polynomial
approximation of degree 4 was used with an approximation error about $5\times
10^{-4}$. Wang et al. \cite{Wang.etal2008} proposed a fast implementation based
on piecewise linear approximation for the QEA of
finite-elongation-nonlinear-elastic (FENE) model
of flexible polymer. Recently, a fast evaluation algorithm for the Bingham
moments with numerical error less than $5\times 10^{-8}$ was given by Luo et al.
\cite{Luo.etal2017}, where series expansions are used for large Bingham
parameter values and a piecewise cubic (Hermite) interpolation is used for
inner region of Bingham parameters.

Only lower order polynomial approximations are used in the above mentioned
numerical methods. When high accuracy is needed, there methods have to use a
huge number of grid points in the parameter space, which leads to
large memory cost. Otherwise, the large numerical error make it hard to prove
the energy dissipation property of the reduced model rigorously, one has to
seek some particular coarse-grain free energy to prove its
dissipation, we refer to
\cite{hu_new_2007}\cite{xu_quasientropy_2020} for this approach.
In this paper, we design efficient high order methods
for Bingham closure approximation on unit circle and unit sphere using
global polynomial approximations and piecewise polynomial approximations, which
can reduce the implementation error to $10^{-15}$ with much smaller memory
cost. We
hope that with the new
efficient and accurate implementation, the Bingham distribution and related QEA
can be applied to
wider applications including by not confined to the closure approximations of
polymer dynamics.

The rest of this paper is organized as follows. In Section 2, we give a brief
introduction to quasi-equilibrium closure approximation with focus on
antipodally symmetric functions on $n$-sphere.
We then consider the closure approximation on unit circle in Section 3
and consider the closure approximation on unit sphere in Section 4.
A summary with a short discussion on the extension to higher dimensional cases
 is given in Section 5.

\section{Preliminaries on Quasi-equilibrium closure approximation}

\subsection{The moment closure problem}
We take the Doi-Smoluchowski equation that describes the dynamics of rod-like
polymers as
an example to introduce the moment closure problem. For rod-like polymers whose
molecules can be described by an orientation (unit) vector $\m \in
\mathbb{R}^3, |\m|=1$. We use configuration
distribution function $f(\x,\m,t)$ to denote the number density
of polymer molecules located at spatial
position $\x$ with orientation $\m$ at time $t$.
The corresponding dynamics is described by following Doi-Smoluchowski equation
\cite{doi_theory_1986}:
\begin{equation} \label{eq:fp}
\frac{\d f}{\d t}=\frac{1}{De}\mathcal{R}\cdot (\mathcal{R}f+f\mathcal{R}U)-\mathcal{R}\cdot(\m \times \mathcal{\kappa}\cdot \m f),
\end{equation}
where $\mathcal{R}=\m\times \frac{\partial}{\partial\m}$ is the gradient
operator on spherical surface, $De$ is the Deborah number, and
  $\kappa$ is  the (fluid) velocity gradient tensor.
$U$ is the molecular interaction potential, usually taken as Maire-Saupe
potential \cite{maier_einfache_1958}:
\begin{equation}
U(\m, t) = U_{0}\int |\m \times \m^{\prime}|^{2} f(\m^{\prime}, t)\d
\m^{\prime}= U_{0}(1-M):\m\m.
\end{equation}
Here $M=\langle\m\m\rangle$ is the second order moment tensor, $U_0$ is a
constant.
We use shorthand notations `$\cdot$' and `$:$' for tensor contractions, with
the former is an
extension of  inner product of two vectors. More precisely, suppose $A =
(A_{i_1, \ldots,
i_m})$, $B = (B_{i_1, \ldots, i_n})$, then $A\cdot B = \big(\sum_{k} A_{i_1,
\ldots, i_{m-1}, k} B_{k, j_{2}, \ldots, j_{n}}\big) $ and $A : B =
\big(\sum_{k_1,  k_2} A_{i_1,
\ldots, i_{m-2}, k_1, k_2} B_{k_1, k_2, j_{3}, \ldots, j_{n}}\big) $.
Putting multiple $\m$ next to each other means tensor product, e.g. $\m\m =
(m_i m_j)$, where $\m=(m_i)$.
Note that for simplicity,
the spatial variation terms are ignored in Eq. \eqref{eq:fp}. We refer to
\cite{yu_kinetic_2007} and \cite{yu_nonhomogeneous_2010a}
 for a
detailed description of the
nonhomogeneous model and the effects of anisotropic spatial diffusion.

Multiplying Eq.(\ref{eq:fp}) by $\m\m$ and then integrating both sides of the
resulting equation with respect to $\m$ on unit sphere, we obtain the evolution
equation for second-order moment tensor $M$, which involves fourth-order
moment tensor $Q=\langle\m\m\m\m\rangle$:
\begin{align}\label{eq:ev_fp}
\begin{split}
\frac{\d M}{\d t} & = \frac{1}{De}\int_{|\m|=1}\m\m \mathcal{R}\cdot(\mathcal{R}f+f\mathcal{R}U)\d \m
-\int_{|\m|=1}\m\m \mathcal{R}\cdot(\m\times\mathcal{\kappa}\cdot\m f)\d\m\\
& = \frac{1}{De}[\langle \mathcal{R}\cdot\mathcal{R}(\m\m)\rangle+U_{0}
\langle\mathcal{R}(\m\m)\cdot \mathcal{R}(\m\m)\rangle:M]
-\mathcal{\kappa}:\langle\m\m\times\mathcal{R}(\m\m)\rangle\\
& = \frac{1}{De}\big[-6(M-\frac{I}{3})+4U_{0}(M\cdot
M-M:Q)\big]+\mathcal{\kappa}\cdot
M +M\cdot \mathcal{\kappa}^{T}-\mathcal{\kappa}:Q.
\end{split}
\end{align}
The system (\ref{eq:ev_fp}) is not closed due to the existence of
higher order moments $Q$, which are also unknown.
To close the system, we must approximate $Q$ using the seconder-order moments
$M$. This process is known as \emph{ closure approximation}.

\subsection{The Quasi-equilibrium approximation}
The QEA uses a distribution that maximizes entropy with observed
information as constraints to close the dynamical system.
Given some lower order moments, the maximum entropy approximation of a
distribution is formulated as (see e.g. \cite{mead_maximum_1984}):
\begin{equation} \label{eq:max-entropy}
\max_f S[f], \quad  S[f] := -\int_{\m\in \Omega}  \big(f(\m)\ln f(\m) - f(\m)\big)
\mathrm{d} \m,
\end{equation}
subject to
\begin{equation}\label{eq:moments}
\int_{\m\in \Omega} p_j(\m)f(\m)\d \m = P_j, \quad \text{for}\  j=0,1,\cdots,k,
\end{equation}
where $\m$ is the microscopic configuration variable defined in
 $\Omega\in \mathbb{R}^n$. $f(\m)$ denotes a configurational density
function on $\Omega$.
Note that in \cite{mead_maximum_1984}, $p_j(\m)$ are assumed to be monomials,
here $\{p_j(\m)\}$ are linearly independent polynomials of $\m$, which allows
the
use of orthogonal polynomials to get better \emph{numerical stability} for
$k>1$.
$P_j$ is the moment of $f$ corresponding to $p_j$ for every $j$. Here we take
$p_0(\m)=1$, $P_0=1$ due to the normalization condition of the CDF. The system
\eqref{eq:max-entropy}-\eqref{eq:moments}, which is a well-posed concave
maximization problem, can be solved by the Lagrange
multiplier method. Define the Lagrangian as
\begin{equation} \label{eq:LagF}
L[f;\lambda]=S[f]+\sum_{j=0}^{k}\lambda_j\big(\int p_j(\m)f(\m)\d \m - P_j\big),
\end{equation}
where $\{\lambda_j \}_{j=0}^{k}$ are the Lagrange multipliers. Then the
solution
to the problem \eqref{eq:max-entropy}-\eqref{eq:moments} satisfies
\begin{equation}
\frac{\delta L}{\delta f} = \ln f + \sum_{j=0}^{k}\lambda_j p_j(\m)=0,
\end{equation}
from which we obtain
\begin{equation}\label{eq:MED}
f(\m)=f_{\lambda}(\m):=\exp\bigg(-\sum_{j=0}^{k}\lambda_j p_j(\m)\bigg)
=\frac{1}{z}\exp\bigg(-\sum_{j=1}^{k}\lambda_j p_j(\m)\bigg).
\end{equation}
Here the normalization constant $z$ is a function of $\{\lambda_j \}_{j=1}^{k}$
defined by
\begin{equation}\label{eq:partitionF}
z=z(\lambda_1, \cdots,
\lambda_k):=\int_{\Omega}\exp\Big(-\sum_{j=1}^{k}\lambda_{j}
p_j(\m)\Big)\mathrm{d} \m.
\end{equation}
The PDF given in form \eqref{eq:MED} is known as a maximum entropy
distribution(MED) or quasi-equilibrium distribution(QED).
One crucial property of such a QEA is that it keeps the free energy dissipation
law of the original dynamical system, see e.g.
\cite{ilg_canonical_2002}\cite{yu_nonhomogeneous_2010a}.

The Lagrange
multipliers $\{\lambda_j \}_{j=0}^{k}$ are determined by Eq.
\eqref{eq:moments} and \eqref{eq:MED}.
The function $z(\lambda_1, \cdots, \lambda_k)$ defined in \eqref{eq:partitionF}
is known as the partition function. It carries all the information of the
distribution function $f_\lambda(\m)$. For example, by taking derivative of
$z$ with respect to $\lambda_j$, we obtain
\begin{equation}
\frac{\partial z}{\partial
\lambda_j}=-\int_{\Omega}\exp\bigg(-\sum_{j=1}^{k}\lambda_j
p_j(\m)\bigg)p_j(\m)\mathrm{d} \m=-zP_j,
\end{equation}
i.e.
\begin{equation}
P_j=-\frac{1}{z}\frac{\partial z}{\partial \lambda_j}=-\frac{\partial
\ln(z)}{\partial \lambda_j},
\quad j=1, \cdots, k.
\end{equation}

To apply maximum entropy distribution \eqref{eq:MED} to close Eq.
\eqref{eq:ev_fp}, one needs to first find $\{\lambda_j \}_{j=0}^{k}$ for given
$\{P_j \}_{j=1}^{k}$ by solving \eqref{eq:moments} and \eqref{eq:MED} together,
then evaluate $Q$ by its definition
\begin{equation}\label{eq:Qexpress}
    Q_j=\int_{\Omega} f_{\lambda}(\m)q_j(\m) \d\m, \quad j=1,\ldots, n_Q.
\end{equation}

It is obvious that solving \eqref{eq:moments} and \eqref{eq:MED} to find the
inverse mapping from the lower-order moments to Lagrange multipliers and
evaluating the integration are computationally expensive. Fortunately,
the mapping between the given moments and the Lagrange multipliers are
smooth functions, we may pre-calculate the integrations at some grid
points and use
interpolation to fast evaluate the inverse mapping at other points.
Note that, one may use the dual approach, which takes Lagrange multipliers as
variables and derive corresponding evolution equations for them from the
Fokker-Planck equation. In both approaches, the closure approximations are not
avoidable. Since the moments of CDF usually have special
physical meaning and are physically measurable, we use in this paper the
standard approach which uses lower order moments as evolution variables.

Note that, even though the maps from lower order moments to high order moments
are quite smooth in QEA, the computational cost grows very fast for large $n$
and $k$. So, in this paper we will consider only numerical implementations for
the cases with $n=2,3$ and the given information is second order moment
tensor, and
leave the cases with larger values of $n$
and $k$ for a future study.

\subsection{Basic mathematical properties of QEA for symmetric
distributions}

We present here some basic theoretical results. We first introduce
some definitions.
\begin{definition}[Antipodally symmetric domain]
    A domain $\Omega \in \mathbb{R}^n$ is said to be \emph{antipodally
    symmetric} if
    for any $\m \in \Omega$, then $-\m \in \Omega$.
\end{definition}
\begin{definition}[Antipodally symmetric function/distribution]
    A function/distribution defined on an antipodally symmetric domain $\Omega
    \in \mathbb{R}^n$ is said to be \emph{antipodally symmetric} if
    $f(\m) =f(-\m)$, for all $\m \in \Omega$.
\end{definition}
\begin{definition}[canonical domain]
    A domain $\Omega \in \mathbb{R}^n$ is said to be \emph{canonical} if
    for any $\m \in \Omega$, $U \m \in \Omega$, where $U$ is an arbitrarily
    given orthogonal matrix.
\end{definition}

\begin{remark}
Note that, the entire Euclid space $\mathbb{R}^n$, the unit
       circle, sphere, hyperspheres, $n$-dimensional ball and spherical annulus
       are all
       canonical. There are more geometries that are antipodally symmetric,
       e.g.
       $n$-dimensional hypercube, ellipsoidal surface/ball, etc.
\end{remark}

For antipodally symmetric distributions, we have the following observation.
\begin{lemma}\label{lm:1}
	Suppose $f(\m)$ is an antipodally symmetric function defined on $\Omega \in
	\mathbb{R}^n$, and $p(\m)$ is a
    monomial of total degree $k$, where $k$ is odd. Then we have
	\begin{equation}\label{eq:oddmeq}
	\int_{\Omega}f(\m)p(\m)\d \m=0.
	\end{equation}
\end{lemma}
\begin{proof}
	By the definition of antipodally symmetric function, we have
	\begin{equation}
	\int_{\Omega}f(\m)p(\m)\d\m = \int_{\Omega}f(-\m)p(-\m)\d\m=(-1)^k \int_{\Omega}f(\m)p(\m)\d\m.
	\end{equation}
	which leads to Eq.\eqref{eq:oddmeq}. \qed
\end{proof}

Lemma \ref{lm:1} says that for antipodally symmetric
distributions, all the odd order moments are zero. Since the zeroth order
moment is a
normalization constant, the first set of nonzero moments are the second order
moments. If the polynomials $\{p_{j} \}_{j=1}^{k}$ are the monomials of
degree 2 in Eq. \eqref{eq:moments} and \eqref{eq:MED}, then we can rewrite the
equations as
\begin{equation}\label{eq:moment2nd}
M=\int_{\Omega}f_{B}(\m)\m \m \d\m,
\end{equation}
\begin{equation}\label{eq:B2nd}
f_{B}(\m)=\frac{1}{z(B)}\exp(-\m^tB\m).
\end{equation}
where $B$ and $M$ are two $n$-dimensional symmetric second order tensors.
\begin{lemma}\label{thm:2.2}
	Let matrix $M$ and $B$ satisfy \eqref{eq:moment2nd}-\eqref{eq:B2nd}.
	If $\Omega \in \mathbb{R}^n$ is canonical,
	then $M$ and $B$ are diagonalizable simultaneously.
\end{lemma}
\begin{proof}
	Note that $M$ and $B$ are both symmetric matrices, so there exist an
	orthogonal matrix $U$ and a diagonal matrix $\Lambda $, such that
	$B=U^t\Lambda U$, i.e. $\Lambda=UBU^t$. With Eq.\eqref{eq:moment2nd},
	\begin{equation}
	M = \int_{\Omega} \frac{1}{z}\exp(-\m^t U^t \Lambda U \m)\m\m \d\m
	= \int_{\Omega} \frac{1}{z}\exp(-(\m')^t\Lambda \m')U^t\m' \m' U \d\m'
	\end{equation}
	where $\m' = U\m \in \Omega$. Thus, by rewriting $\m'$ as $\m$, we have
	\begin{equation}\label{eq:diagM}
	UM U^t = \int_{\Omega} \frac{1}{z}\exp(-\m^t\Lambda \m)\m \m \d\m.
	\end{equation}

        For the
	off-diagonal elements in $\m\m$, such as $m_{i}m_{j}$, $i \neq j$, we
	make a variable change $v_{k}=m_{k} $, for $k\neq j$, and $v_{j}=-m_{j}$.
	Since ${\Lambda}$ is diagonal, we have
	\begin{equation}
	\frac{1}{z}\int_{\Omega}\exp(-\m^t {\Lambda}\m)m_i m_j\d\m=
	-\frac{1}{z}\int_{\Omega}\exp(-{\bf v }^t {\Lambda}{\bf v} ) v_i v_j\d{\bf
	v},
	\end{equation}
	from which we obtain that the off-diagonal elements are zero, i.e. $\Sigma
	= U M U^t $ is a diagonal matrix. Therefore, $M$ and $B$ can be
	diagonalized simultaneously. \qed
\end{proof}

\begin{remark}
If $\Omega$ is a tensor-product domain, e.g. $\Omega = \mathbb{R}^n$, then
according to Lemma \ref{thm:2.2}, after the diagonalization, the $n$
components of $\m$ are decoupled to $n$ 1-dimensional problems, which makes
the corresponding moment closure approximation an easy task.
However, if $\Omega$ is not of tensor-product type, such as the unit
$n$-sphere, we can't obtain decoupled sub-problems. We will focus on the latter
case where $\Omega$ is $n$-sphere and the corresponding canonical antipodally
symmetric distribution $f_B$ defined in \eqref{eq:B2nd} are known as the Bingham
distribution\cite{bingham_antipodally_1974}. We will consider
the numerical implementation of the cases $n=2$ and $3$ in next two sections.
\end{remark}

\section{Symmetric QEA on unit circle}
In this section, we consider the moment closure problem \eqref{eq:moment2nd}
and
\eqref{eq:B2nd}, where $\Omega=\{\m \in \mathbb{R}^2: |\m|=1\}$. Given second
order moments $M$, we study the moment closure problem that how to fast
calculate
\eqref{eq:Qexpress}, where $Q=\{ Q_j, j=1, \ldots, n_Q \}$ are forth order
moments.

\subsection{Some theoretical results}
By Lemma \ref{thm:2.2}, $M$ and $B$ can be diagonalized simultaneously.
So, we consider the diagonalized case first.
Notice that the distribution function $f_{B+cI}(\m)$ is identical to
$f_{B}(\m)$, where the shift $cI$ with $I$ being identity matrix goes into
the normalization constant $z$.
Hence, we need only consider the case where $B={\diag}(\la,-\la)$
with
$ \lambda \geq 0 $. Using polar coordinates $\m=(m_1,
m_2)=(\cos\theta,\sin\theta)$, where
$0 \leq \theta \leq 2\pi $, we have
$f_{B}(\m)=f_{\lambda}(\m):=\frac{1}{z}\exp(\la\cos(2\theta))$, and
\begin{equation}\label{eq:mij}
m_{ij}=\frac{1}{z(\lambda)}\int_{0}^{2\pi}
\exp(\la\cos(2\theta))\m_{i}\m_{j}\mathrm{d}\theta, \quad i,j=1,2,
\end{equation}
where
\begin{equation}
z(\lambda)=\int_{0}^{2\pi}
\exp(\la\cos(2\theta))\d\theta, \label{eq:defz}
\end{equation}
and $\{m_{ij}\}_{i,j=1,2}$ are the elements of second-order moment $M$.
By the definition of $\Omega$ and Lemma \ref{lm:1}, we have $m_{ij}=0$ if
$i\neq j$ and
$m_{11}+m_{22}=1$. So we only
need one free variable to define the second order moments. We take this
variable as
\begin{equation} \label{eq:defmu}
\mu(\lambda):=m_{11}-m_{22}=\frac{1}{z(\lambda)}\int_{0}^{2\pi}
\exp(\la\cos(2\theta))\cos(2\theta)\d\theta.
\end{equation}
It is obvious that $z(0)=2\pi$, $z(+\infty)=+\infty$, and $\mu(\lambda)$ can be
represented by $z(\lambda)$ as
\begin{equation}\label{eq:3.3}
\mu(\lambda)=\frac{z^{\prime}(\lambda)}{z(\lambda)}.
\end{equation}
Furthermore
\begin{equation}\label{eq:defmuprime}
\mu^\prime(\lambda)=\frac{z^{\prime\prime}(\lambda)}{z(\lambda)}
-\big(\frac{z^\prime(\lambda)}{z(\lambda)}\big)^2
=\frac{z^{\prime\prime}(\lambda)}{z(\lambda)}-\mu^2(\lambda).
\end{equation}
The elements of the fourth order moments are defined as
\begin{equation}\label{eq:qijkl}
q_{ijkl}=\frac{1}{z(\lambda)}\int_{0}^{2\pi}
\exp(\la\cos(2\theta))\m_{i}\m_{j}\m_{k}\m_{l}\d\theta, \quad i,j,k,l=1,2,
\end{equation}
which can be classified as
\begin{equation}\label{eq:qm}
q_{m}=\frac{1}{z(\lambda)}\int_{0}^{2\pi}
\exp(\la\cos(2\theta))\cos^{m}(\theta)\sin^{4-m}(\theta)\mathrm{d}\theta,
\quad m=0,1,2,3,4.
\end{equation}
It is easy to verify that $q_{1}=q_{3}=0$,
$q_{0}+q_{2}=\frac{1-\mu}{2}$, $q_{2}+q_{4}=\frac{1+\mu}{2}$.
So there is only one independent variable. We take it as
\begin{equation}\label{eq:defeta}
\eta(\lambda)=\frac{z^{\prime\prime}(\lambda)}{z(\lambda)}=\frac{1}{z(\lambda)}\int_{0}^{2\pi}
\exp(\la\cos(2\theta))\big(\cos(2\theta)\big)^{2}\mathrm{d}\theta.
\end{equation}
It follows from Eq.\eqref{eq:defmuprime} and Eq.\eqref{eq:defeta} that
\begin{equation}\label{eq:eta_mu}
\mu^{\prime}(\lambda)=\eta(\lambda)-\mu^{2}(\lambda).
\end{equation}
We will use this relation in Newton's method in next subsection.

If we know the values of $\mu$ and $\eta$, then we can calculate the elements
in the second order and forth order moments according to the following equations:
\begin{align}
\label{eq:re_m_mu}
& m_{12}  =  m_{21} = 0, \quad m_{11} =  \frac{1+\mu}{2},  \quad  m_{22} =
\frac{1-\mu}{2},\\
\label{eq:re_q_eta}
& q_{1}  =  q_{3} = 0, \quad q_{0}  = \frac{1-2\mu+\eta}{4}, \quad q_{2}  =
\frac{1-\eta}{4},  \quad q_{4} = \frac{1+2\mu+\eta}{4}.
\end{align}

If we define
\begin{equation}\label{eq:defzk}
z^{(k)}(\lambda):=\frac{\mathrm{d}z^{(k-1)}(\lambda)}{\mathrm{d}\lambda}
=\int_{0}^{2\pi}\exp\big({\lambda}\cos(2\theta)\big){\cos^k(2\theta)}\mathrm{d}\theta,\quad
 k=1,2,\cdots,
\end{equation}
then we can easily obtain the following conclusion.
\begin{lemma}
	The partition function $z(\lambda)$ is an analytic function with all
	derivatives uniformly bounded. More precisely, we have for finite $\lambda$
	\begin{equation*}z(\lambda)>z^{(2)}(\lambda)>z^{(4)}(\lambda)
    >\cdots>z^{(2k)}(\lambda) \rightarrow 0^{+},\quad \mbox{for}\ k\rightarrow
    \infty,
	\end{equation*}
	and $0<z^{(2k+1)}(\lambda)<z^{(2k)}(\lambda)$ for any integer $k\ge 0$.
\end{lemma}

\begin{theorem}\label{th:3.1}
	The $\mu^\prime(\lambda)$ in Eq.\eqref{eq:defmuprime} is positive for
	any finite $\lambda$.
\end{theorem}
\begin{proof}
	By using Eqs.\eqref{eq:defz},
	\eqref{eq:defzk}, and the Cauchy-Schwarz inequality, we have
	\begin{align*}
	(z^{\prime}(\lambda))^2 & = \big(\int_{0}^{2\pi}\exp(\la\cos(2\theta))\cos(2\theta)\d\theta\big)^2\\
	& \leq \int_{0}^{2\pi}\exp(\la\cos(2\theta))\d\theta \cdot \int_{0}^{2\pi}\exp(\la\cos(2\theta))(\cos(2\theta))^2\d\theta\\
	& = z(\lambda)z^{\prime \prime}(\lambda)
	\end{align*}
	The equal sign hold if and only if $z(\lambda)=z^{\prime \prime}(\lambda)$, which can't be true for finite $\lambda$. So we obtain that
	\begin{equation*}
	(z^{\prime}(\lambda))^2 < z(\lambda)z^{\prime \prime}(\lambda),
	\end{equation*}
	i.e.
	\begin{equation}
	(\mu(\lambda))^2 < \eta(\lambda),
	\end{equation}
	which is equivalent to $\mu^\prime(\lambda)>	 0$ by Eq.
	\eqref{eq:defmuprime}.
 \qed
\end{proof}

\subsection{Evaluation of $z(\lambda)$, $\mu(\lambda)$ and $\eta(\lambda)$ as
functions of $\lambda$}

From Eq.(\ref{eq:defz}), we have
that \begin{equation}
z(\lambda)=\int_{0}^{2\pi}\exp(\la\cos2\theta)\d\theta=2\pi I_0(\la),
\end{equation}
where $I_0(\lambda)$
is the modified Bessel function of the first kind. $I_p(\lambda)$ has an
integral representation \cite[page 376]{abramowitz_handbook_1972}
\begin{equation*}
I_p(\lambda)=\frac{1}{\pi}\int_{0}^{\pi}\exp(\lambda\cos\theta)\cos(p\theta)\d\theta.
\end{equation*}
Similarly, from Eq. (\ref{eq:defmu})(\ref{eq:defeta}), we have
\begin{equation}\label{eq:mulambda}
\mu(\lambda)=\frac{1}{z(\lambda)}\int_{0}^{2\pi}\exp(\la\cos2\theta)\cos2\theta\d\theta
=\frac{I_1(\la)}{I_0(\la)},
\end{equation}
\begin{equation}\label{eq:etalambda}
\eta(\lambda)=\frac{1}{z(\lambda)}\int_{0}^{2\pi}
\exp(\la\cos(2\theta))\frac{\cos(4\theta)+1}{2}\mathrm{d}\theta=\frac{1}{2}
[\frac{I_2(\la)}{I_0(\la)}+1].
\end{equation}
So to calculate $z(\lambda), \mu({\lambda}), \mu({\lambda})$ efficiently, we
only need a fast subroutine to calculate the modified Bessel function of the
first kind, which can be done by using series expansion, differential equation,
or continued fractions \cite{olver_nist_2010}. It is implemented in
most mathematical software and  libraries, e.g. {\tt Netlib},
Python's  {\tt{SciPy}} library, Matlab etc.
To obtain high accurate numerical results, we adopt the multiple-precision
implementation $besseli(p,\lambda)$ in {\tt
MPmath}\cite{johansson_mpmath_2020}, which  is a Python library
for real and complex floating-point arithmetic with arbitrary precision.

\subsection{Representation of 4-th order moments in terms of 2-nd order moments}

Now we describe how to represent all 4-th order moments in terms of 2-nd order
moments. This involves several steps.
\begin{enumerate}
    \item The first step is to diagonalize the
    second order moment tensor by using an orthogonal transform, then to
    calculate the variable $\mu$ from second order moments by Eq.
    \eqref{eq:defmu}, i.e. $\mu= m_{11} - m_{22}$.
    \item In the second step, we approximate
    $\eta(\mu)$ by an truncated Legendre polynomial approximation, sine $\eta$
    as a function of $\mu$ is very smooth.
    With Legendre coefficients pre-calculated, this evaluation can be done very
    efficiently and have spectral accuracy.
    \item In the third step, we calculate the forth order
    moments in the transformed coordinates by Eq. \eqref{eq:re_q_eta}, followed
    by an coordinate transform to calculate the 4-th order moments in
    original coordinates.
\end{enumerate}

For the first step, let's denote the second order moment tensor before
diagonalization by
\begin{equation*}
    \hat{M} = \left(\begin{array}{cc}
        a & b\\
        b & 1-a
\end{array}\right),
\end{equation*}
where $a\in [0,1]$ and $\Delta = a(1-a)-b^2 \ge 0$ by the symmetric
semi-positive
definite property. Then the two eigenvalues are given by
\begin{equation}
    \lambda_{1,2} = \frac{1\pm \sqrt{1-4\Delta}}{2}
    = \frac{1\pm \sqrt{(2a-1)^2 + 4b^2}}{2}. \label{eq:eigv1d}
\end{equation}
The corresponding orthogonal
transformation matrix is given by
\begin{equation}
U = \begin{cases}
\begin{pmatrix}
\frac{b}{\sqrt{b^2+(\lambda_1-a)^2}}\
\frac{b}{\sqrt{b^2+(\lambda_2-a)^2}} \\
\frac{\lambda_1-a}{\sqrt{b^2+(\lambda_1-a)^2}}\
\frac{\lambda_2-a}{\sqrt{b^2+(\lambda_2-a)^2}}
\end{pmatrix}, &
\mbox{if}\ b\neq 0, \\
\begin{pmatrix}
1 & 0\\
0 & 1\\
\end{pmatrix}, &
\mbox{if}\ b=0, a\ge \frac12, \\
\begin{pmatrix}
0 & 1\\
1 & 0\\
\end{pmatrix}, &
\mbox{if}\ b=0, a < \frac12.
\end{cases}
\label{eq:Utrans}
\end{equation}
The diagonalized second-order moment tensor is given by
\begin{equation}
M = \begin{pmatrix}
m_{11}\ 0 \\
0\
m_{22}
\end{pmatrix}
= \begin{pmatrix}
\lambda_{1}\ \ 0 \\
0\ \
\lambda_{2}
\end{pmatrix}.
\end{equation}
Note that we have $\lambda_1 \ge \lambda_2$.

In the second step, we approximate $\eta(\mu)$, where $\mu=m_{11} -
m_{22}\in[0,1]$, by
an truncated Legendre polynomial approximation:
\begin{equation}\label{eq:3.18}
\eta(\mu)=\tilde{\eta}(x) \approx \sum_{k=0}^{n_l}b_{k}L_{k}(x),
\end{equation}
where $x=2\mu -1\in[-1,1]$, and
\begin{equation}
b_{k} =
\dfrac{2k+1}{2}\sum_{j=0}^{N}\tilde{\eta}(x_{j})L_{k}(x_{j})\omega_{j},
\quad k=0,\ldots, n_l.
\label{eq:3.19}
\end{equation}
Here $\{x_j, \omega_j\}_{j=0}^N, N>n_l$ are the Legendre-Gauss quadrature points
and
weights, they are calculated with high accuracy by using the method described
in \cite[page 99]{shen_spectral_2011a}.
To obtain $b_k$, we need know the values of $\tilde{\eta}(x_j) =
\eta(\mu_j)$, $\mu_j=\frac{x_j+1}{2},\ j=0,\ldots, N$. This can be done by
using Newton's method for Lagrange multiplier $\lambda$, since we know how to
calculate
$\mu(\lambda)$ and $\eta(\lambda)$ by Eq. \eqref{eq:mulambda} and
\eqref{eq:etalambda}.
To start the Newton's method, we need good initials. The procedure
to calculate the Legendre coefficients $\{b_k\}$ is given below.

\begin{figure}[htbp]
    \centering{}
    \includegraphics[width=0.6\textwidth]{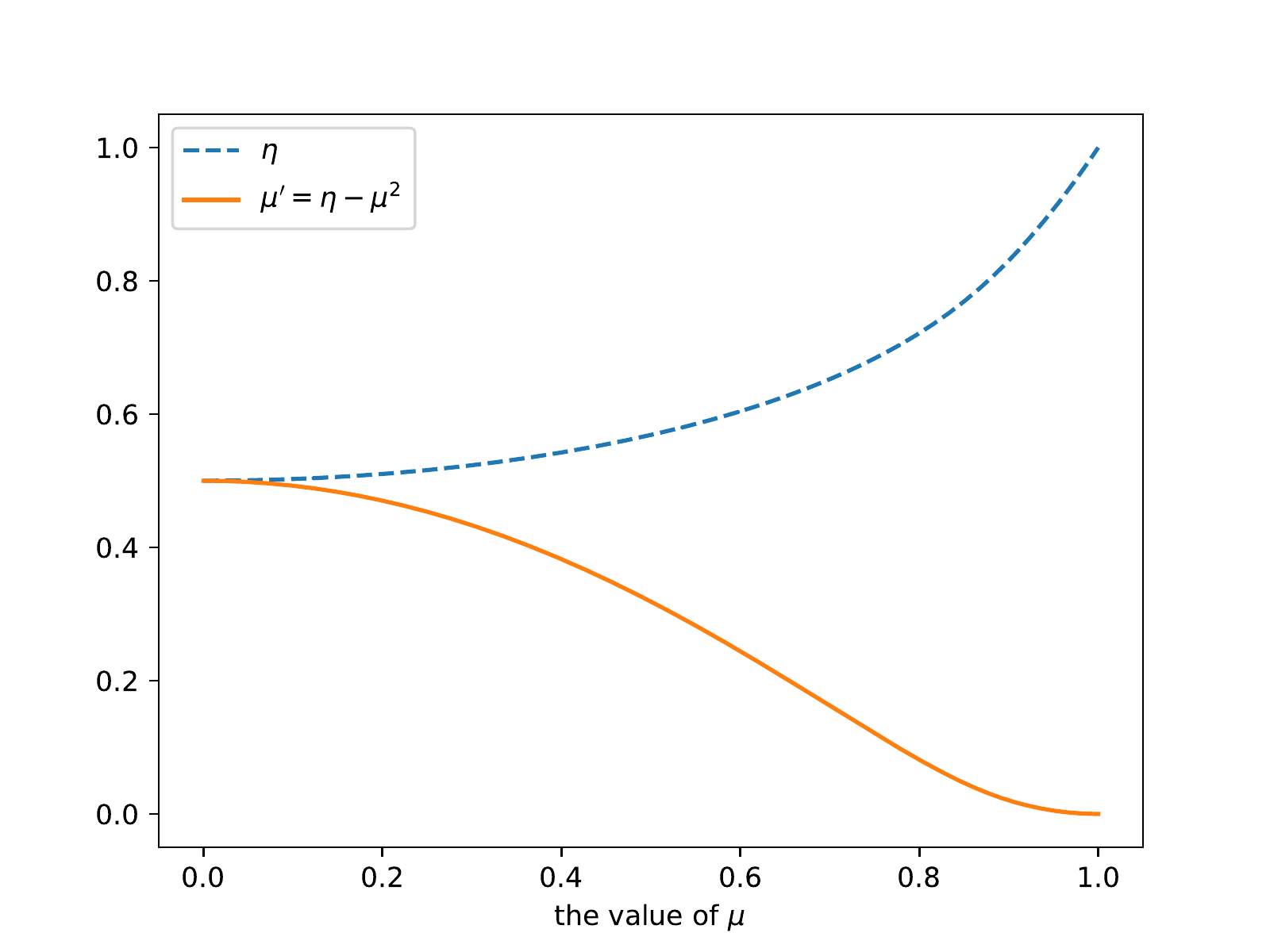}
    \caption{The graphs of $\eta$ and $\eta-\mu^2$ as functions of $\mu$}
    \label{fig_eta_mu}
\end{figure}

\begin{enumerate}[i)]
    \item We first calculate  $\mu(\lambda)$ at a series of $\lambda$ points
    and
    save the results  as $\{\,\lambda_i^0,\mu_i^0\,\},\ i=1,2,\cdots, N$ to
    form a table named $\mu_{init}$.
    \item  For each point $\mu_j$ in $\{\,
    \mu_j\!\!=\!\!\frac{x_j+1}{2}\,\}_{j=0}^N$,
    we find a $\lambda_j^0$ from the table
    $\mu_{init}$ with corresponding $\mu_j^0$ is
    closest to but no larger than $\mu_j$ and use $\lambda_j^0$ as the initial
    value to start
    the following Newton iteration:
\begin{equation}\label{eq:NewtonIter}
\lambda_j^{k+1}=
\lambda_j^k+\dfrac{\mu_j-\mu(\lambda_j^{k})}{\mu^\prime(\lambda_j^k)}=
\lambda_k+\dfrac{\mu_j-\mu(\lambda_j^{k})}{\eta(\lambda_j^{k})
    -\mu^2(\lambda_j^{k})},\quad k=0, 1, \ldots,
\end{equation}
    where relation \eqref{eq:eta_mu} is used to avoid the calculation of
    derivatives.
    We stop the iteration when the distance between
    $\mu(\lambda_j^{k+1})$ and $\mu_{j}$ is no more than $10^{-18}$,
    and the last iteration point $\lambda_j^{k+1}$ is denoted as $\lambda_j^*$.
    \item  For each $j$, calculate $\tilde{\eta}(x_j)=\eta(\lambda_j^*)$ by Eq.
    \eqref{eq:etalambda}.
    Then we can obtain the Legendre expansion coefficients
    $\{b_k\}_{k=0}^{n_l}$ by
    evaluating Eq. \eqref{eq:3.19}.
\end{enumerate}
Note that except for $\mu_j=1$, for which we have $\lambda_j^*=\infty$ and
$\eta(\mu_j)=1$, the Newton iteration  \eqref{eq:NewtonIter} has global
convergence.
To see this, we formally subtract the fixed point $\lambda_j^*$ from both sides
of
\eqref{eq:NewtonIter} to obtain
\begin{equation}
    e_j^{k+1} = e_j^k \left(1 -
    \frac{\mu'(\bar\lambda_j^k)}{\mu'(\lambda_j^k)}\right),
\end{equation}
where $e_j^{k} = \lambda_j^k - \lambda_j^*$ and
$\bar\lambda_j^k\in[\lambda_j^k, \lambda_j^*]$. We have used the mean
value theorem. By Theorem \ref{th:3.1}, $\mu'(\lambda_j^k)$ and
$\mu'(\bar\lambda_j^k)$ are positive. To show the global convergence, we only
need to verify that $\mu'(\lambda)$ is a decreasing function, in such case
$\{\lambda_j^k, j=0, \ldots \}$ is a bounded increasing sequences.
Mathematically proving that $\mu'(\lambda)$ is a decreasing function seems
tedious.
On the other side, we can easily
deduce this result by looking at Figure \ref{fig_eta_mu}, where the graphs of
$\eta$ and $\mu'=\eta - \mu^2$ with respect to $\mu$ are plotted. We see
$\mu'$ as a function of $\mu$ is decreasing, then by the fact that
$\mu(\lambda)$ is an
increasing function, we known  $\mu'$ as a function of $\lambda$ is decreasing.

\begin{figure}[htbp]
    \centering{}
    \includegraphics[width=\textwidth]{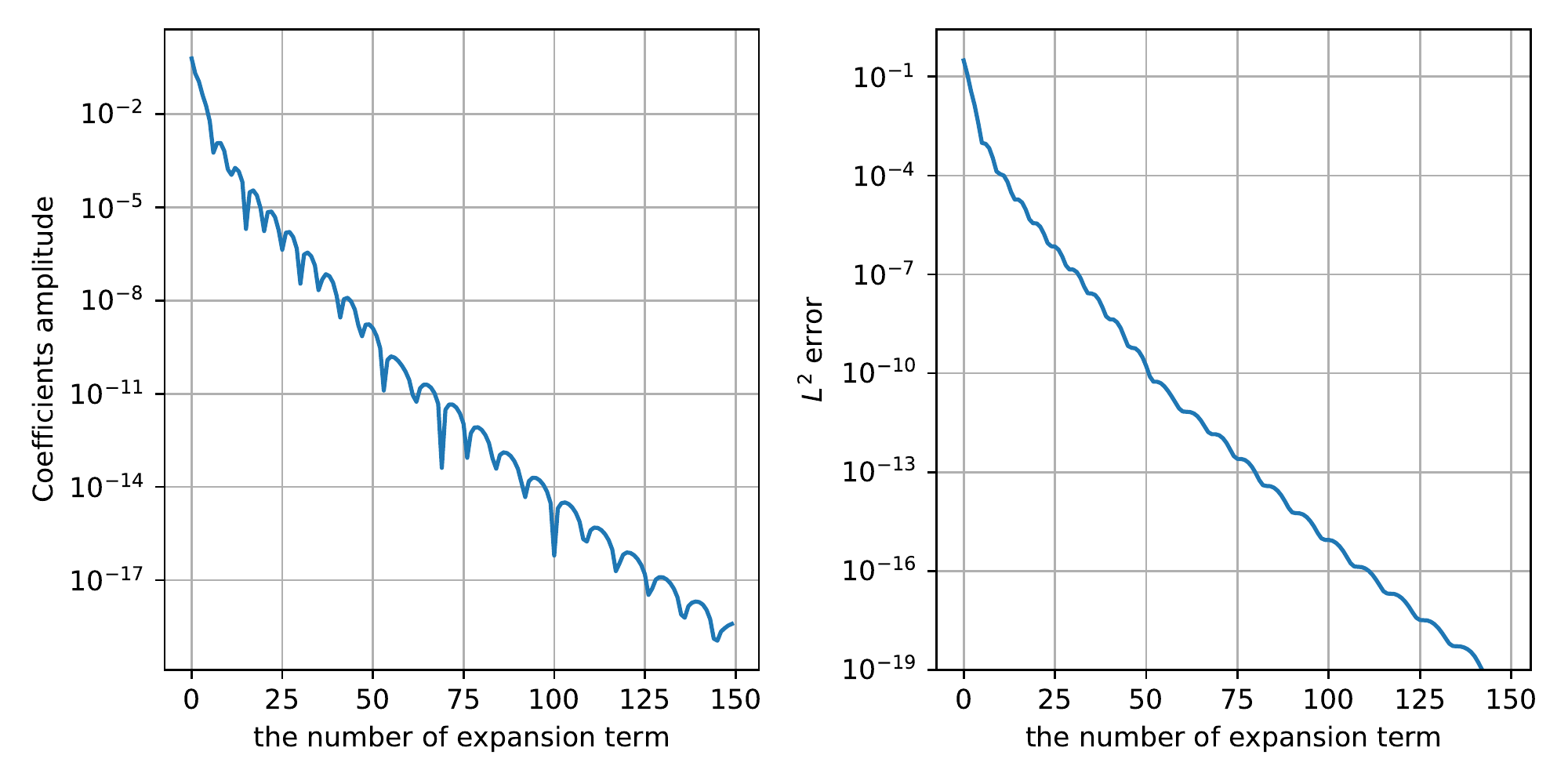}
    \caption{The coefficients amplitude(left) and the estimated $L^2$ error
    (right) of
    Legendre
        expansion of $\eta(\mu)$}\label{fig2}
\end{figure}

In our numerical experiments the Newton iteration usually terminate
in 3 or 4 steps. In Figure \ref{fig2}, we show the amplitudes of Legendre
coefficients calculated using above procedure where we take $n_l=150$ and
$N=200$
in \eqref{eq:3.18} and \eqref{eq:3.19}.
From this figure, we see that the Legendre expansion has spectral
accuracy and the error is about $10^{-18}$ when
125 expansion terms are used.

After the coefficients of the Legendre expansion of $\eta(\mu)$
are pre-calculated, we can easily obtain the elements in the transformed fourth
order moment by Eq. \eqref{eq:re_q_eta}. Then, as the last step of moment
closure approximation, the components of the forth order moments in original
coordinates can be calculated by:
\begin{equation}\label{eq:qhattrans}
\hat{q}_{ijkl} = \sum_{i',j',k',l'}U_{ii'}U_{jj'}U_{kk'}U_{ll'}{q}_{i'j'k'l'},
\end{equation}
where $U_{ij}$ are the components of the transform matrix $U$ defined in Eq.
\eqref{eq:Utrans}.

The overall moment closure procedure for symmetric QEA on unit circle is
summarized in Algorithm \ref{alg:repre_q}.
\begin{algorithm}
	\setstretch{1.35}      
	\caption{Calculate the fourth order moments for symmetric QEA on unit
	circle}
	\label{alg:repre_q}
	{\bf Input:}                      
	The values of the elements in second order moment  $\hat{M}$.\\
	{\bf Output:}                     
	The values of the elements in fourth order moment $\hat{Q}$. 
	\begin{algorithmic}[1]
		\REQUIRE Pre-calculated coefficients $\{b_k, k=0,\ldots, n_l\}$ of
		the Legendre expansion of
		$\eta(\mu)$.
		\STATE Transform $\hat{M}$ into a diagonalized matrix
		$M=U^T\hat{M}U$, where $U$ is given by \eqref{eq:Utrans}.
		\STATE Calculate $\mu=m_{11} - m_{22}$ where $m_{11}$ and $m_{22}$ are
		the diagonal elements in $M$.
		\STATE Map the value of $\mu$ to $x \in [-1,1]$ with $x=2\mu-1$.
		\STATE Use the Legendre expansion given in Eq.(\ref{eq:3.18}) to
		calculate the value of $\eta$.
		\STATE Evaluate the values of $q_0$, $q_2$ and $q_4$ by
		Eq.(\ref{eq:re_q_eta}), then use the definition \eqref{eq:qm} and
		\eqref{eq:qijkl} to calculate $q_{ijkl}$ and use \eqref{eq:qhattrans}
		to calculate $\hat{Q} = (\hat{q}_{ijkl})$.
	\end{algorithmic}
\end{algorithm}

\begin{remark}
The overall storage cost in Algorithm \ref{alg:repre_q} is $O(n_l+m)$, and the
computational time cost
is $O(mn_l)$, where $n_l$ is the number of Legendre coefficients and
$m$ is the number of fourth-order moment tensors to be evaluated.
For $n_l \le 125$, the evaluation of \eqref{eq:3.18} can be done by fast
matrix-vector product on modern computers.
\end{remark}

Note that it is possible to further reduce the computational time cost by
using
a piecewise polynomial approximation. For example, if we divide the range of
$\mu \in$ [0,1] into 6 intervals: $[0,0.5], [0.5,0.73], [0.73,0.84],
[0.84,0.91], [0.91,0.96], [0.96,1]$. Then the numbers of Legendre expansion
terms $n_l$ can be significantly reduced for maintaining similar $L^2$
error of approximating $\eta(\mu)$ in each interval. The result is shown in
Table
\ref{tab:1.1}. We see double precision is achieved on all six intervals
with less than 20 Legendre coefficients, which leads to an overall
computational time cost reduction by a factor of 5. One can divide the interval
into more pieces to further reduce the computational cost. For simplicity, we
will not present more results here.
\renewcommand\arraystretch{1.1}
\begin{table}[ht!]
    \setlength{\tabcolsep}{4mm}{
        \begin{center}
            \caption{Piecewise Legendre approximation of $\eta(\mu)$.}
            \label{tab:1.1}
            \begin{tabular}{lcc}
                \hline
                Interval & $n_l$ & $L^2$ error \\ \hline
                $[0,0.5]$  & 19  & 7.34E-18 \\ \hline
                $[0.5,0.73]$  & 19  & 2.30E-18 \\ \hline
                $[0.73,0.84]$  & 18  & 5.45E-18 \\ \hline
                $[0.84,0.91]$  & 18  & 7.64E-18 \\ \hline
                $[0.91,0.96]$  & 18  & 8.97E-18 \\ \hline
                $[0.96,1]$  & 18  & 5.41E-18 \\ \hline
            \end{tabular}
    \end{center}}
\end{table}

\section{Symmetric QEA on unit sphere}

The overall moment closure procedure based on symmetric QEA on unit sphere is
similar
to that on unit circle. For simplicity we first consider
the diagonalized case.

We consider the diagonalized 3-dimensional problem on unit sphere $\Omega =
\{\,\m \in \mathbb{R}^3,\, |\m|=1\,\}$.
We consider two different cases: the uniaxial case, where two of the three
eigenvalues of $B$ tensor are equal, and the biaxial case where three
eigenvalues of $B$ tensor are all different. Since $f_{B+cI}$ is identical
to $f_B$, we can use a particular shift $cI$ to make the $B$ in uniaxial case
be $\diag(0,0,-2\lambda)$, make the $B$ in biaxial case
$\diag(-\lambda_1-\lambda_2, -\lambda_1+\lambda_2, 0)$.
Here $\lambda_{1}\geq\lambda_{2} \geq 0$.
We first consider the uniaxial case, which is easier to implement.

\subsection{The uniaxial case}
For the uniaxial case,  we set
$B=\mathrm{diag}(0, 0, -2\lambda)$, $\lambda \in \mathbb{R}$. For $\lambda>0$,
$f_B(\m)$ is an oblate
distribution, while $\lambda<0$, it is prolate.
By using spherical coordinates
$\m=(\sin\theta\cos\varphi,\sin\theta\sin\varphi,\cos\theta)$
with $0\leq\theta\leq\pi$, $0\leq\varphi\leq2\pi$, we have
\begin{equation*}
m_{ij}=\int_{0}^{\pi}\int_{0}^{2\pi}\frac{1}{z}
\exp[-2\la\cos^{2}\theta]
\m_{i}\m_{j} \sin\theta\mathrm{d}\theta\mathrm{d}\varphi, \quad i,j=1,2,
\end{equation*}
where
\begin{equation}\label{def_z_5.1}
z(\lambda)=\int_{0}^{\pi}\int_{0}^{2\pi}\exp[-2\la\cos^{2}\theta]
\sin\theta\d\theta\d\varphi.
\end{equation}
It is easy to check that $m_{ij}=0$ if $i\neq j$. The nonzero terms left are
$m_{11}$, $m_{22}$ and
$m_{33}$. We define
\begin{equation}\label{def_mu_5.1}
\mu(\la):=-\frac{z^{\prime}(\lambda)}{z(\lambda)}
=\frac{1}{z(\lambda)}\int_{0}^{\pi}\int_{0}^{2\pi}
\exp[-2\la\cos^{2}\theta]2\cos^2\theta\sin\theta\d\theta\d\varphi.
\end{equation}
Then
\begin{equation}\label{eq:muprime1}
\mu^\prime(\lambda)=-\frac{z^{\prime\prime}(\lambda)}{z(\lambda)}
+\big(\frac{z^\prime(\lambda)}{z(\lambda)}\big)^2
=\mu^2(\lambda)-\frac{z^{\prime\prime}(\lambda)}{z(\lambda)}.
\end{equation}
And the second order moments are related to $\mu$ by
\begin{equation}
\label{eq:m_mu}
m_{33} = \frac{\mu}{2},  \quad m_{11} = m_{22} = \frac{2-\mu}{4}.
\end{equation}

The fourth order moments are defined as
\begin{equation*}
q_{ijkl}=\int_{0}^{\pi}\int_{0}^{2\pi}\frac{1}{z(\lambda)}
\exp[-2\la\cos^{2}\theta]
\m_{i}\m_{j}\m_{k}\m_{l}\sin\theta\d\theta\d\varphi,\quad i,j,k,l=1,2,3.
\end{equation*}
One may check that only the
following several terms: $q_{1111}$, $q_{2222}$, $q_{3333}$, $q_{1122}$,
$q_{1133}$, $q_{2233}$ are nonzero, and they satisfy following constraints
\begin{align}\label{eq:relmq}
\left(\begin{array}{ccc}
1 & 1 & 0\\
1 & 0 & 1\\
0 & 1 & 1
\end{array}\right)\left(\begin{array}{c}
q_{1122}\\
q_{1133}\\
q_{2233}
\end{array}\right)=\left(\begin{array}{c}
m_{11}-q_{1111}\\
m_{22}-q_{2222}\\
m_{33}-q_{3333}
\end{array}\right).
\end{align}
Actually, there is only one independent variable in forth order moments. Let's
define it as
\begin{equation}\label{def_eta_5.1}
\eta(\lambda)=\frac{z^{\prime\prime}(\lambda)}{z(\lambda)}
=\frac{1}{z(\lambda)}\int_{0}^{\pi}\int_{0}^{2\pi}
\exp[-2\la\cos^{2}\theta]
4\cos^4\theta \sin\theta \d\theta\d\varphi.
\end{equation}
It follows from \eqref{eq:muprime1} that
\begin{equation}\label{eq:eta_mu_5.1}
\mu^{\prime}(\lambda) =\mu^{2}(\lambda) - \eta(\lambda).
\end{equation}
Similar to Theorem \ref{th:3.1}, one can prove that $\mu^{\prime}(\lambda)$
is positive by using the Cauchy-Schwarz inequality.
By using the definitions, the symmetry between $m_1$ and $m_2$, the
relation \eqref{eq:relmq}, we find that fourth order moments are related to
$\eta$ and $\mu$ by:
\begin{align}
\label{eq:q_eta}
& q_{3333} = \frac{\eta}{4},
\quad q_{1133}  = q_{2233} = 
\frac{\mu}{4}-\frac{\eta}{8},
\quad q_{1111} = q_{2222} = 3 q_{1122} = \frac{3}{8}(1-\mu)+\frac{3}{32}\eta.
\end{align}

To efficiently evaluate $z(\lambda)$, we rewrite it by using Eq.
(\ref{def_z_5.1}) as
\begin{align}
z(\lambda) &
=4\pi \int_{0}^{\frac{\pi}{2}}
\exp[-2\la\cos^{2}\theta]\sin\theta\d\theta
\nonumber\\
& =4\pi\int_{0}^{1}
\exp[-2\la t^{2}] \d t  \nonumber\\
&
=2\pi\int_{0}^{1}
\exp[-2\la x] x^{-\frac{1}{2}}\d
x \nonumber\\
& =2\pi \frac{\Gamma(1)\Gamma(\frac{1}{2})}{\Gamma(\frac{3}{2})}
 \ _1F_1(\frac12;\frac{3}{2};-2\la)
\nonumber\\
& =4\pi
\ _1F_1(\frac12;\frac{3}{2};-2\la), \label{eq:zf11}
\end{align}
where
\begin{equation}\label{eq:chf11}
_1F_1(a;b;\la)=\frac{\Gamma(b)}{\Gamma(a)\Gamma(b-a)}
\int_{0}^{1}\exp(\la x)x^{a-1}(1-x)^{b-a-1}\d x
\end{equation}
is the confluent hypergeometric function \cite[Chapter 13]{olver_nist_2010}.
Similarly, we have
\begin{equation}
\mu(\lambda)=\frac{8\pi}{3z}
\ _1F_1(\frac32;\frac{5}{2};-2\la),
\end{equation}
\begin{equation}
\eta(\lambda)=\frac{16\pi}{5z}
\ _1F_1(\frac52;\frac{7}{2};-2\la).
\end{equation}
We use the function {\tt hyp1f1(a,b,z)} in {\tt MPmath}
\cite{johansson_mpmath_2020}
to calculate the confluent hypergeometric function $_1F_1$
with high accuracy.

\begin{figure}[htbp]
	\centering{}
	\includegraphics[width=\textwidth]{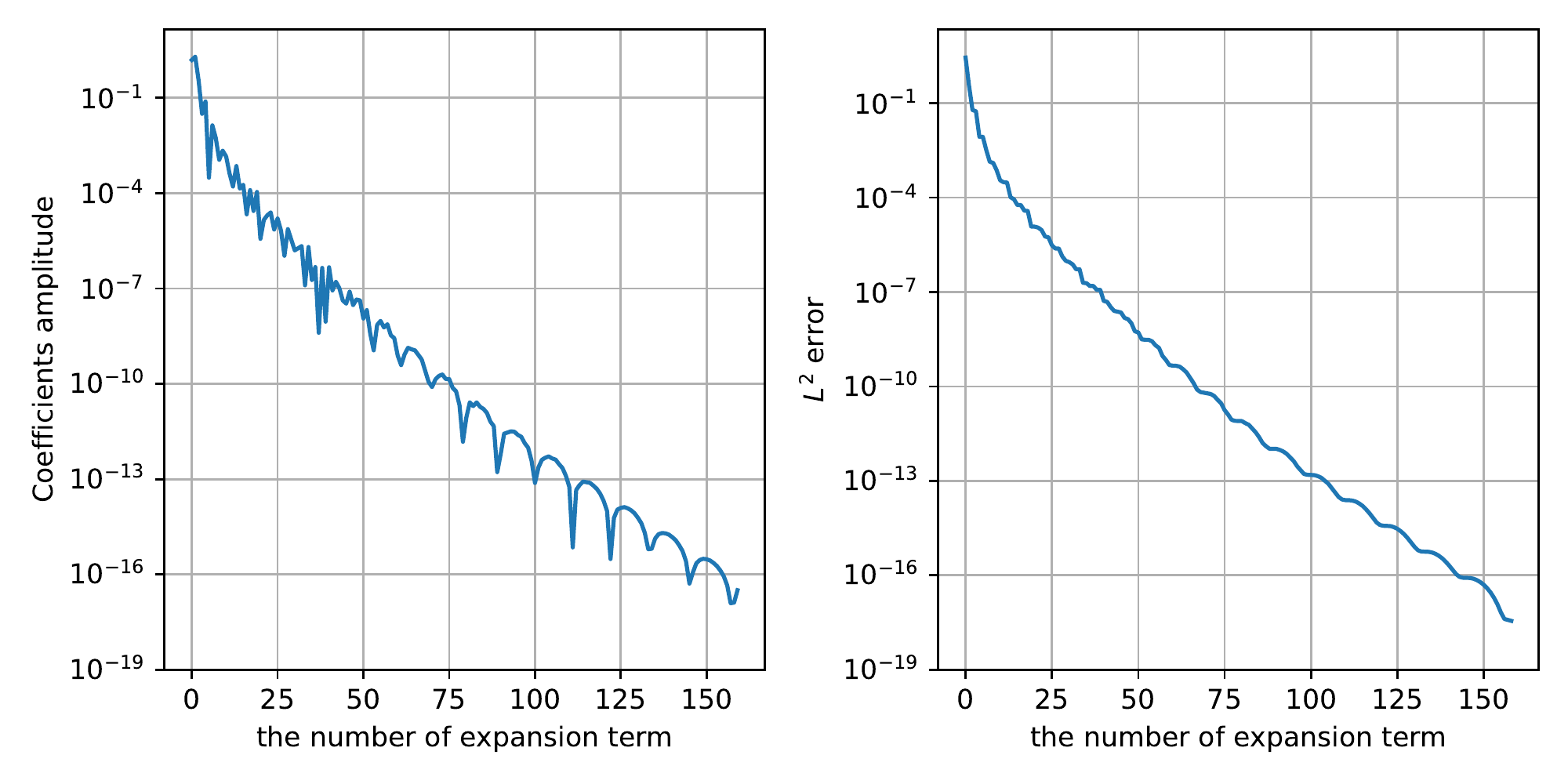}
	\caption{The coefficients (left) and the $L^2$ error (right) of Legendre expansion of $\eta$}\label{fig5.1}
\end{figure}

Then, as done in the 2-dimensional case,
we consider an Legendre polynomial approximation of $\eta(\mu)$
as a function of $\mu$, which is similar to
Eq.\eqref{eq:3.18},  but here $x=\mu-1$ since $\mu \in [0, {2}]$.
We use again Newton's method to find the $\lambda$ values
that produce $\mu(\lambda)$ on mapped Legendre-Gauss points, then use these
$\lambda$
values to calculate corresponding $\mu$ and $\eta$ values, and use them to
obtain Legendre coefficients.
We take $n_l=160$ and use $N=200$ Legendre-Gauss points to
compute the coefficients of the Legendre polynomial approximation.
The result is given in Figure \ref{fig5.1}, from which we see that the
Legendre expansion has spectral
accuracy and the $L^2$ error is reduced to about $10^{-16}$ (smaller than
double precision $2.2\times 10^{-16}$) with
less than 150 expansion terms are used.

After the coefficients of the Legendre approximation of $\eta(\mu)$ is
obtained, we can efficiently calculate $\eta$ for given $\mu$, and then
obtain the elements in fourth order moment by \eqref{eq:q_eta}.
The overall procedure is very similar to Algorithm
\ref{alg:repre_q} for the two dimensional case.

About 150 global Legendre terms are needed to obtain double precision, which
means the computational cost for evaluating $m$ fourth order moments is about
$O(n_l m)$, where $n_l\approx 150$. We can further reduce the computational
cost by
using piecewise polynomial approximations.
For example, we can divide the range of $\mu \in$ $[0, {2}]$ into
multiple intervals,
then use Legendre expansions
to approximate $\eta(\mu)$ on each interval. The number of expansion terms
$n$ and the $L^2$  approximation errors on each interval are showed in Table
\ref{table:2.1}. 
We see that with only 18 expansion terms, the $L^2$ approximation
error can be reduced to less than $4\times 10^{-17}$ on all intervals.
Note that $\mu=2/3$ corresponds to the special case where $m_{11} = m_{22} = m_{33}$.
So the first 5 intervals in Table \ref{table:2.1} are oblate cases, while the last 6 intervals
are prolate cases.

\renewcommand\arraystretch{1.5}
\begin{table}
	\setlength{\tabcolsep}{7mm}{
	\caption{Piecewise Legendre approximation for the uniaxial Bingham
	closure on
	sphere.}
    \label{table:2.1}
    \begin{center}
        \begin{tabular}{lcc}
		\hline
		Interval & $n_l$ & $L^2$ error \\ \hline
		$[0, 0.045]$  & 18  & 1.41E-17 \\ \hline
		$[0.045, 0.103]$  & 18  & 2.34E-17 \\ \hline
		$[0.103, 0.2]$  & 18  & 2.88E-17 \\ \hline
		$[0.2,0.36]$  & 18  & 1.69E-17 \\ \hline
		$[0.36,\frac{2}{3}]$  & 18  & 3.42E-17 \\ \hline
		$[\frac{2}{3},1.26]$  & 18  & 1.28E-17 \\ \hline
		$[1.26,1.56]$  & 18  & 1.13E-17 \\ \hline
		$[1.56, 1.73]$  & 18  & 1.26E-17 \\ \hline
		$[1.73, 1.84]$  & 18  & 1.38E-17 \\ \hline
		$[1.84,1.925]$  & 18  & 1.80E-17 \\ \hline
		$[1.925,2]$  & 18  & 1.89E-17 \\ \hline
	\end{tabular}
    \end{center}
    }
\end{table}

\subsection{The biaxial case}

In biaxial case, we take
$B = \diag(-\lambda_1-\lambda_2, -\lambda_1+\lambda_2, 0)$, where
$\lambda_{1}>\lambda_{2} > 0$. Note that the limit cases
$\lambda_2=0$ and $\lambda_2 \rightarrow \infty$ are reduced to uniaxial
distribution.
By using spherical coordinates
$\m=(\sin\theta\cos\varphi,\sin\theta\sin\varphi,\cos\theta)$
with $0\leq\theta\leq\pi$, $0\leq\varphi\leq2\pi$, we have
\begin{equation*}\label{eq:m3d}
    m_{ij}=\int_{0}^{\pi}\int_{0}^{2\pi}
    \frac{1}{z}\exp[-(\lambda_{1}+\lambda_{2}\cos(2\varphi))\sin^{2}\theta]
    \m_{i}\m_{j} \sin\theta\mathrm{d}\theta\mathrm{d}\varphi, \quad i,j=1,2,
\end{equation*}
where
\begin{equation}\label{3d_def_z}
z(\lambda_{1},\lambda_{2})=\int_{0}^{\pi}\int_{0}^{2\pi}\exp[-(\lambda_{1}+\lambda_{2}\cos(2\varphi))\sin^{2}\theta]
\sin\theta\mathrm{d}\theta\mathrm{d}\varphi.
\end{equation}
By the definition of $B$, we have $m_{11} \le m_{22} \le m_{33}$.
It is easy to check that $m_{ij}=0$ if $i\neq j$ by symmetry
or direct integration. The nonzero second order moments are $m_{11}$,
$m_{22}$ and $m_{33}$, with constraint
$m_{11}+m_{22}+m_{33}=1$, so we have two independent variables.
We take them as
\begin{align}
  \mu_{1}(\lambda_{1},\lambda_{2})&=-\frac{\partial z}{\partial\lambda_{1}}\frac{1}{z} =\frac{1}{z}\int_{0}^{\pi}\int_{0}^{2\pi}\exp[-(\lambda_{1}+\lambda_{2}\cos(2\varphi))\sin^{2}\theta]\sin^{3}\theta
                                    \d\theta\d\varphi,\\
  \mu_{2}(\lambda_{1},\lambda_{2})&=-\frac{\partial z}{\partial\lambda_{2}}\frac{1}{z} =\frac{1}{z}\int_{0}^{\pi}\int_{0}^{2\pi}\exp[-(\lambda_{1}+\lambda_{2}\cos(2\varphi))\sin^{2}\theta]
                                    \cos(2\varphi)\sin^{3}\theta\d\theta\d\varphi.
\end{align}
Variables $\mu_1, \mu_2$ are related to second order moments $m_{11}, m_{22}$ by
\begin{equation}\label{eq:m2mu}
    \mu_1 = m_{11}+m_{22},\quad
    \mu_2 = m_{11}-m_{22}.
\end{equation}
The fourth order moments are defined as
\begin{equation*}
q_{ijkl}=\int_{0}^{\pi}\int_{0}^{2\pi}\frac{1}{z(\lambda)}
\exp[-(\lambda_{1}+\lambda_{2}\cos(2\varphi))\sin^{2}\theta]
\m_{i}\m_{j}\m_{k}\m_{l}\sin\theta\text{d}\theta\text{d}\varphi,\quad
i,j,k,l=1,2,3.
\end{equation*}
It is easy to check that the nonzero terms are: $q_{1111}$, $q_{2222}$,
$q_{3333}$, $q_{1122}$, $q_{1133}$, $q_{2233}$.
They satisfy the relation
\begin{align}\label{eq4.4}
\left(\begin{array}{ccc}
1 & 1 & 0\\
1 & 0 & 1\\
0 & 1 & 1
\end{array}\right)\left(\begin{array}{c}
q_{1122}\\
q_{1133}\\
q_{2233}
\end{array}\right)=\left(\begin{array}{c}
m_{11}-q_{1111}\\
m_{22}-q_{2222}\\
m_{33}-q_{3333}
\end{array}\right).
\end{align}
So there are three independent variables. We define them as
\begin{align}
\eta_{1}(\lambda_{1},\lambda_{2})&=\frac{1}{z}\frac{\partial^{2}z}{\partial\lambda_{1}^{2}}
=\frac{1}{z}\int_{0}^{\pi}\int_{0}^{2\pi}\exp[-(\lambda_{1}+\lambda_{2}\cos2\varphi)\sin^{2}\theta]\sin^{5}\theta
\d\theta\d\varphi\\
\eta_{2}(\lambda_{1},\lambda_{2})&=\frac{1}{z}\frac{\partial^{2}z}{\partial\lambda_{1}\partial\lambda_{2}}
=\frac{1}{z}\int_{0}^{\pi}\int_{0}^{2\pi}\exp[-(\lambda_{1}+\lambda_{2}\cos2\varphi)\sin^{2}\theta]
\cos(2\varphi)\sin^{5}\theta\d\theta\d\varphi\\
\eta_{3}(\lambda_{1},\lambda_{2}) & =\frac{1}{z}\frac{\partial^{2}z}{\partial\lambda_{2}^{2}}
=\frac{1}{z}\int_{0}^{\pi}\int_{0}^{2\pi}\exp[-(\lambda_{1}+\lambda_{2}\cos2\varphi)\sin^{2}\theta]
\cos^{2}(2\varphi)\sin^{5}\theta\d\theta\d\varphi
\end{align}
From the relation between $\mu(\lambda)$ and $\eta(\lambda)$, we derive that
\begin{align}\label{eq5.8}
\begin{split}
\left\{ \begin{array}{cll}
\eta_{1} = & \mu_{1}^{2}-\frac{\partial\mu_{1}}{\partial\lambda_{1}}, &\\
\eta_{2}  = & \mu_{1}\mu_{2}-\frac{\partial\mu_{1}}{\partial\lambda_{2}}
=& \mu_{1}\mu_{2}-\frac{\partial\mu_{2}}{\partial\lambda_{1}},\\
\eta_{3}  = & \mu_{2}^{2}-\frac{\partial\mu_{2}}{\partial\lambda_{2}},&
\end{array}\right.
\end{split}
\end{align}
which will be used to compute the Jacobi matrix in the Newton's method.

\begin{theorem}\label{th:4.1}
	The Jacobi matrix of $\frac{\partial(\mu_1, \mu_2)}{\partial(\lambda_1,
	\lambda_2)}$ is negative semi-definite.
\end{theorem}
\begin{proof}

    To show the Jacobi matrix $\frac{\partial(\mu_1, \mu_2)}{\partial(\lambda_1,
        \lambda_2)}$ is negative semi-definite, we first define
    function $f(\la_1,\la_2)=\log(z(\la_1,\la_2))$. Since
    \begin{align*}
    \left(\begin{array}{cc}
    \frac{\partial^{2}f(\la_1,\la_2)}{\partial\lambda_{1}^{2}} &
    \frac{\partial^{2}f(\la_1,\la_2)}{\partial\lambda_{1}\partial\lambda_{2}}\\
    \frac{\partial^{2}f(\la_1,\la_2)}{\partial\lambda_{1}\partial\lambda_{2}} &
    \frac{\partial^{2}f(\la_1,\la_2)}{\partial\lambda_{2}^{2}}
    \end{array}\right)
    =-\left(\begin{array}{cc}
    \frac{\partial\mu_{1}}{\partial\lambda_{1}} &
    \frac{\partial\mu_{1}}{\partial\lambda_{2}}\\
    \frac{\partial\mu_{2}}{\partial\lambda_{1}} &
    \frac{\partial\mu_{2}}{\partial\lambda_{2}}
    \end{array}\right),
    \end{align*}
    we only need to show that $f(\lambda_1, \lambda_2)$ is a convex function,
    or to show that for any given $\lambda_1\ge \lambda_2 \ge 0$, $a^2+b^2=1$,
    $f(\lambda_1+\gamma a, \lambda_2 + \gamma b)$ as a function of $\gamma$ is
    convex. By direct calculation, we have
    \begin{equation*}
        \frac{\d ^2 f}{\d \gamma^2} = \frac{1}{z^2}
        \Big[ z \frac{\d ^2 z}{\d \gamma^2} -
        (\frac{\d z}{\d \gamma})^2\Big],
    \end{equation*}
    where
    \begin{equation}\label{z_gamma}
    z(\gamma)=\int_{0}^{\pi}\int_{0}^{2\pi}\exp[
    -(\lambda_{1}+\lambda_{2}\cos(2\varphi)
    + \gamma (a+b\cos 2\varphi))
    \sin^{2}\theta]
    \sin\theta\mathrm{d}\theta\mathrm{d}\varphi.
    \end{equation}
    Similar to Theorem \ref{th:3.1}, by using Cauchy-Schwartz inequality,
    we have $(\d z/ \d \gamma)^2 \le z(\d^2 z/\d \gamma^2)$, which means
    $f$ as a function of $\gamma$ is convex. The theorem is proved.
   \qed

\end{proof}

Given the values of $\mu_1$, $\mu_2$ and $\eta_1$, $\eta_2$, $\eta_3$,
the second order and forth order moments of the biaxial Bingham distribution
can be obtained by
\begin{align}
\label{eq:re_mandmu}
& m_{11} = \frac{\mu_{1}+\mu_{2}}{2}, \quad m_{22}=\frac{\mu_{1}-\mu_{2}}{2},
\quad  m_{33} =1-\mu_{1},\\
\begin{split}
\label{eq:re_qandeta}
& q_{1111}  =  \frac{\eta_{1}+2\eta_{2}+\eta_{3}}{4}, \quad q_{1122}  = \frac{\eta_{1}-\eta_{3}}{4}, \quad q_{2222}=\frac{\eta_{1}-2\eta_{2}+\eta_{3}}{4},\\
& q_{1133}  = \frac{(\mu_{1}+\mu_{2})-(\eta_{1}+\eta_{2})}{2},  \quad q_{2233}=\frac{(\mu_{1}-\mu_{2})-(\eta_{1}-\eta_{2})}{2},\\
& q_{3333} = 1-2\mu_{1}+\eta_{1}.
\end{split}
\end{align}

Now we describe how to efficiently calculate $z, \mu_1, \mu_2, \eta_1, \eta_2,
\eta_3$ and do the moment closure approximation.

Similar to the uniaxial case, the partition function and moments can be written
as integrations of confluent hypergeometric functions.
Since for large and close $\lambda_1,
\lambda_2$ values, the integrands $_1F_1(a;b;
-(\lambda_1+\lambda_2\cos(2\varphi)))$ are localized
at $\cos(2\varphi)\approx -1$,
we use Legendre-Gauss quadrature to do numerical integration
in $\varphi$
variable to put more grid points near $\cos(2\varphi)\approx -1$.
To this end, we write those quantities as:
\begin{equation}\label{re_z}
z(\la_1, \la_2)=2\int_{0}^{2\pi}\!
_1F_1(1;\frac{3}{2};-(\lambda_{1}+\lambda_{2}\cos2\varphi))\d\varphi
= 2\pi\int^{1}_{-1}\!
_{1}F_{1}(1;\frac{3}{2};-(\lambda_{1}+\lambda_{2}\cos\pi t))\d t
\end{equation}

\begin{align}
& \mu_{1}(\la_{1},\la_{2})
=\frac{4\pi}{3z}\int_{-1}^{1} \
_1F_1(2;\frac{5}{2};-(\lambda_{1}+\lambda_{2}\cos\pi t))\d t,\\
& \mu_{2}(\la_{1},\la_{2})
=\frac{4\pi}{3z}\int_{-1}^{1} \ _1F_1(2;\frac{5}{2};
-(\lambda_{1}+\lambda_{2}\cos\pi t))\cos({\pi}t)\d t,\\
& \eta_{1}(\la_{1},\la_{2})
=\frac{16\pi}{15z}\int_{-1}^{1} \ _1F_1(3;\frac{7}{2};
-(\lambda_{1}+\lambda_{2}\cos\pi t))\d t,\\
& \eta_{2}(\la_{1},\la_{2})
=\frac{16\pi}{15z}\int_{-1}^{1} \ _1F_1(3;\frac{7}{2};
-(\lambda_{1}+\lambda_{2}\cos\pi t))\cos({\pi}t)\d t,\\
& \eta_{3}(\la_{1},\la_{2})
=\frac{16\pi}{15z}\int_{-1}^{1} \ _1F_1(3;\frac{7}{2};
-(\lambda_{1}+\lambda_{2}\cos\pi t))
{\cos^2(\pi t)}\d t,
\end{align}
where $_1F_1$ is defined in \eqref{eq:chf11}, they are evaluated using
the function {\tt hyp1f1(a,b,z)} in {\tt
MPmath} to get high accuracy. Since the integrands are all even functions,
we can use half the Gauss points to save computational time.

\begin{figure}[htbp]
	\centering
	\begin{subfigure}[The region of $\mu_1$ and $\mu_2$]
		{\begin{tikzpicture}[scale=6]
		\shade[top color=gray!100,bottom color=blue]
		(0,0) -- (1/2,-1/2) -- (2/3,0);
		\fill (0.5,-0.5) circle(0.3pt) node[below]{$(\frac{1}{2},-\frac{1}{2})$};

		\draw[->] (-0.1,0) -- (0.75,0)  node[right]{$\mu_1$};
		\draw[->] (0,-0.6) -- (0,0.1) node[above] {$\mu_2$};
		\draw[step=0.1, dashed] ( 0,-0.6 ) grid ( 0.7, 0);

		\foreach \x/\xtext in {0.5/\frac{1}{2}, 0.6667/\frac{2}{3}}
		\draw[shift={(\x,0)}] (0pt,-0.8pt) -- (0pt,0.3pt) node[above] {$\xtext$};

		\foreach \y/\ytext in {-0.5/-\frac{1}{2}}
		\draw[shift={(0,\y)}] (0.8pt,0pt) -- (-0.3pt,0pt) node[left] {$\ytext$};

		\fill (0,0) circle(0.4pt) node[above=6pt,left]{$A$};
		\fill (0.5,-0.5) circle(0.4pt) node[above=6pt,right]{$B$};
		\fill (2/3,0) circle(0.4pt) node[above=6pt,right]{$C$};

		\end{tikzpicture}}
	\end{subfigure}
	\hspace{2.5em}
	\begin{subfigure}[The region of $x$ and $y$]
		{\begin{tikzpicture}[scale=1.8]
		\shade[top color=gray!50,bottom color=blue!100]
		(-1,-1) -- (1,-1) -- (1,1) -- (-1,1);

		\draw[->] (-1.2,0) -- (1.2,0) node[right] {$x$};
		\draw[->] (0,-1.2) -- (0,1.2) node[above] {$y$};

		\foreach \x/\xtext in {-1/-1, 1/1}
		\draw[shift={(\x,0)}] (0pt,1pt) -- (0pt,-1pt) node[right=4pt,below] {$\xtext$};

		\foreach \y/\ytext in {-1/-1,  1/1}
		\draw[shift={(0,\y)}] (1pt,0pt) -- (-1pt,0pt) node[left=3pt,below] {$\ytext$};

		\fill (-1,-1) circle(1.2pt) node[below=6pt,left]{$A$};
		\fill (1,-1) circle(1.2pt) node[below=6pt,right]{$B$};
		\fill (-1,1) circle(1.2pt) node[above=6pt,left]{$C$};
		\fill (1,1) circle(1.2pt) node[above=6pt,right]{$C^{\prime}$};

		\end{tikzpicture}}
    \end{subfigure}
	\caption{The transformation between the triangular and rectangular domains}
	\label{fig:5.1}
\end{figure}
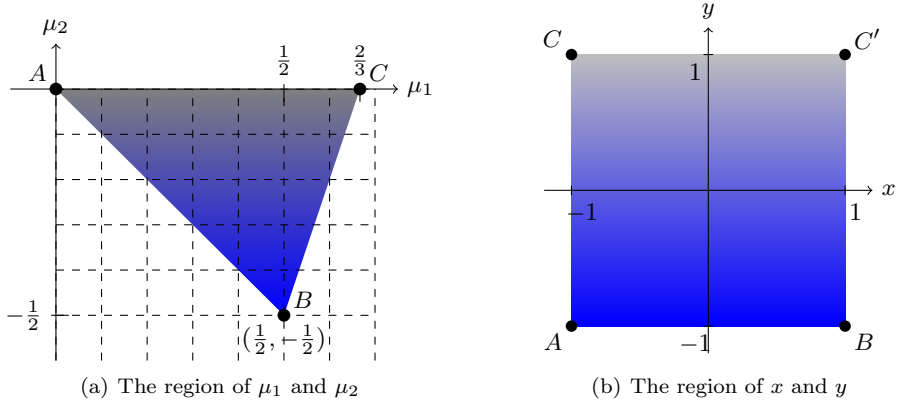

We have shown how to calculate $\mu_1, \mu_2$, $\eta_1, \eta_2$ and $\eta_3$
for given $\lambda$. Next,
we show how to efficiently calculate $\eta_1, \eta_2$ and $\eta_3$ for given
$\mu_1, \mu_2$.
Similar to the 2-dimensional case, we use Legendre method to approximate
functions $\tilde{\eta}_i(\mu_1, \mu_2), i=1,2,3$.
Since the domain of $(\mu_1, \mu_2)$, shown in Figure \ref{fig:5.1}(a),
is not of tensor product form, we first transform it into standard domain
$[-1,1]^2$ by using following mapping
\begin{align}\label{eq:map2d}
\left\{ \begin{array}{ccc}
\mu_{1}  &= & \frac{(1+x)(1-y)}{8}+\frac{y}{3}+\frac{1}{3},\\
\mu_{2}  &= & -\frac{(1+x)(1-y)}{8}.
\end{array}\right.
\end{align}
The corresponding inverse mapping is
\begin{align}\label{tr_mu_x}
\left\{ \begin{array}{ccc}
x & = & \frac{8\mu_{2}}{3\mu_{1}+3\mu_{2}-2}-1,\\
y & = & 3\mu_{1}+3\mu_{2}-1.
\end{array}\right.
\end{align}
where $x\in[-1,1]$, $y\in[-1,1]$.
Then we use Legendre-Gauss points in the transformed domain to
calculate the Legendre approximation coefficients. Denote the Legendre
approximation of $\tilde{\eta}_i(\mu_1, \mu_2)$ as
\begin{align}\label{eq5.11}
\eta_i(\lambda_1, \lambda_2) = \tilde{\eta}_i(\mu_{1},\mu_{2})=\hat{\eta}(x,y)
 \approx \sum_{s=0}^{n_1}\sum_{t=0}^{n_2}b^i_{st}L_{s}(x)L_{t}(y),
 \quad i=1,2,3,
\end{align}
where
\begin{align}\label{eq5.12}
\begin{split}
b^i_{st}
&=\frac{1}{\gamma_{s}\gamma_{t}}\sum_{i=0}^{N_1}\sum_{j=0}^{N_2}\hat{\eta}(x_{i},y_{j})
L_{s}(x_{i})L_{t}(y_{j})\omega_{i} \omega_j,
\end{split}
\end{align}
with $\gamma_{s}=\frac{2}{2s+1}$, $\gamma_{t}=\frac{2}{2t+1}$,
$\{x_i, \omega_i\}_{i=0}^{N_1}$ and $\{y_j, \omega_j \}_{j=0}^{N_2}$
are the Legendre-Gauss quadrature points and weights.

Again, we use Newton's method to  obtain the values of $\lambda_1,
\lambda_2$ that
produce Legendre-Gauss points in $x,y$ domain:
\begin{equation}
\lambda^{k+1}=\lambda^{k}-J_{F}(\lambda^{k})^{-1}F(\lambda^{k}) \label{Niter}
\end{equation}
where $\lambda^{k}=(\lambda_{1}^{k},\lambda_{2}^{k})$,
$F(\lambda^{k})=\mu^{*}-\mu(\lambda^{k})$, $\mu^*$ is the image of
Legendre-Gauss point $(x_i, y_j)$ under mapping \eqref{eq:map2d}.
$J_{F}(\lambda^{k})$ is the Jacobi matrix:
\begin{equation*}
J_{F}(\lambda^{k})=-J_{\mu}(\lambda^{k})=-\left.\left(\begin{array}{cc}
\frac{\partial\mu_{1}}{\partial\lambda_{1}} &
\frac{\partial\mu_{1}}{\partial\lambda_{2}}\\
\frac{\partial\mu_{2}}{\partial\lambda_{1}} &
\frac{\partial\mu_{2}}{\partial\lambda_{2}}
\end{array}\right) \right|_{(\lambda_1,\lambda_2)=(\lambda_1^k, \lambda_2^2)}.
\end{equation*}
Then equation (\ref{Niter}) can be rewritten as
\begin{equation*}
\lambda^{k+1}=\lambda^{k}+J_{\mu}(\lambda^{k})^{-1}(\mu^{*}-\mu(\lambda^{k})).
\end{equation*}
The derivatives in Jacobi matrix $J_{\mu}(\lambda)$ can be removed by using
\eqref{eq5.8}. Similar to the 2-dimensional case, we use a table to find
closed $\lambda$ points to initialize Newton's method and the iteration
is terminated if the $L^2$ distance of the objective functions between
two adjacent iterations is smaller than a given tolerance. The tolerance we
used here is $10^{-15}$. According to Theorem \ref{th:4.1}, the Newton
iteration is well-defined, the system corresponds to a convex optimization
problem. We expect a global convergence as in the 2-dimensional case.
Our numerical results show that the iterations usually terminate in less than
20 steps.

\begin{figure}[htbp]
	\centering
	\includegraphics[width=\linewidth]{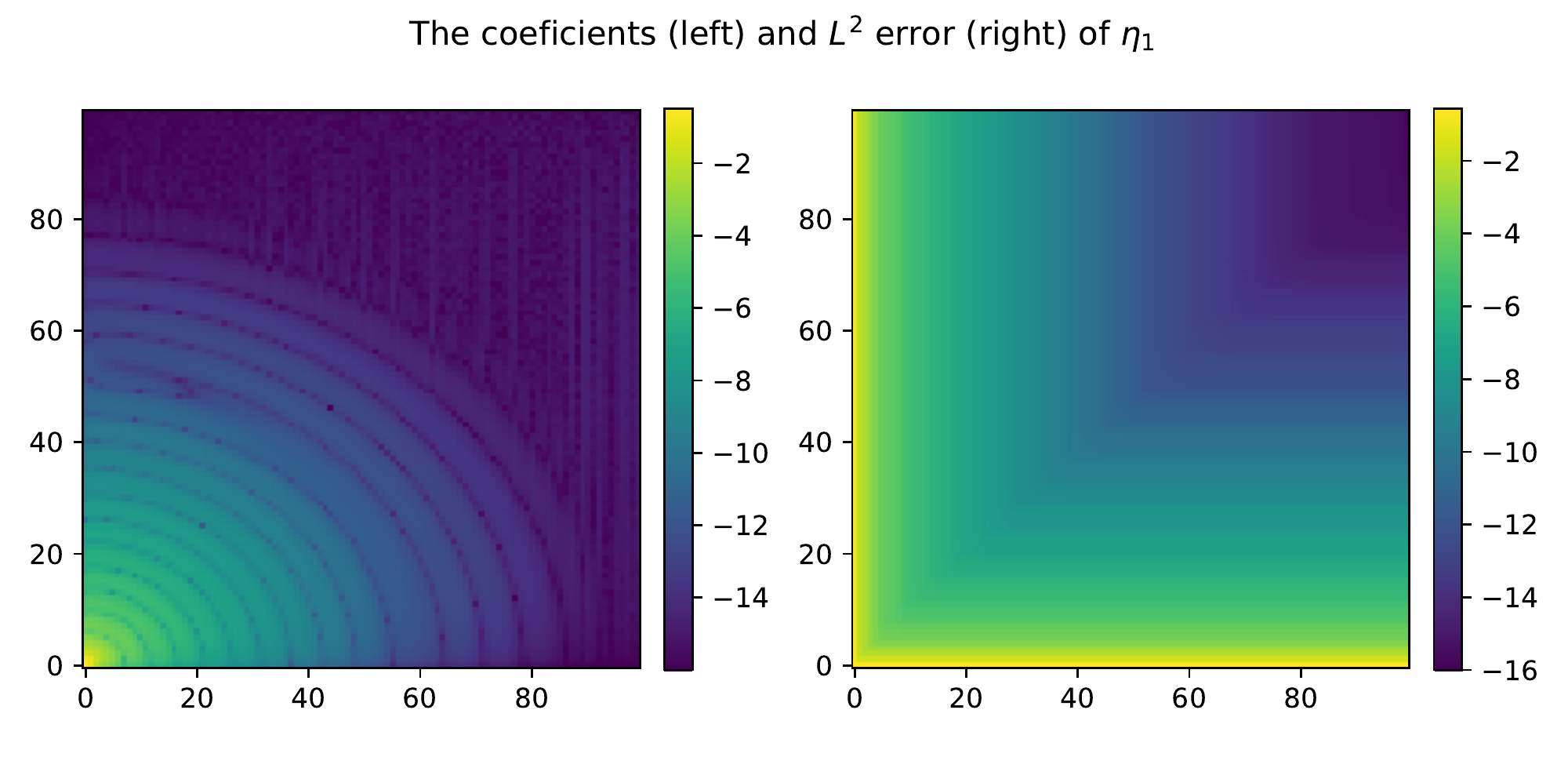}
	\includegraphics[width=\linewidth]{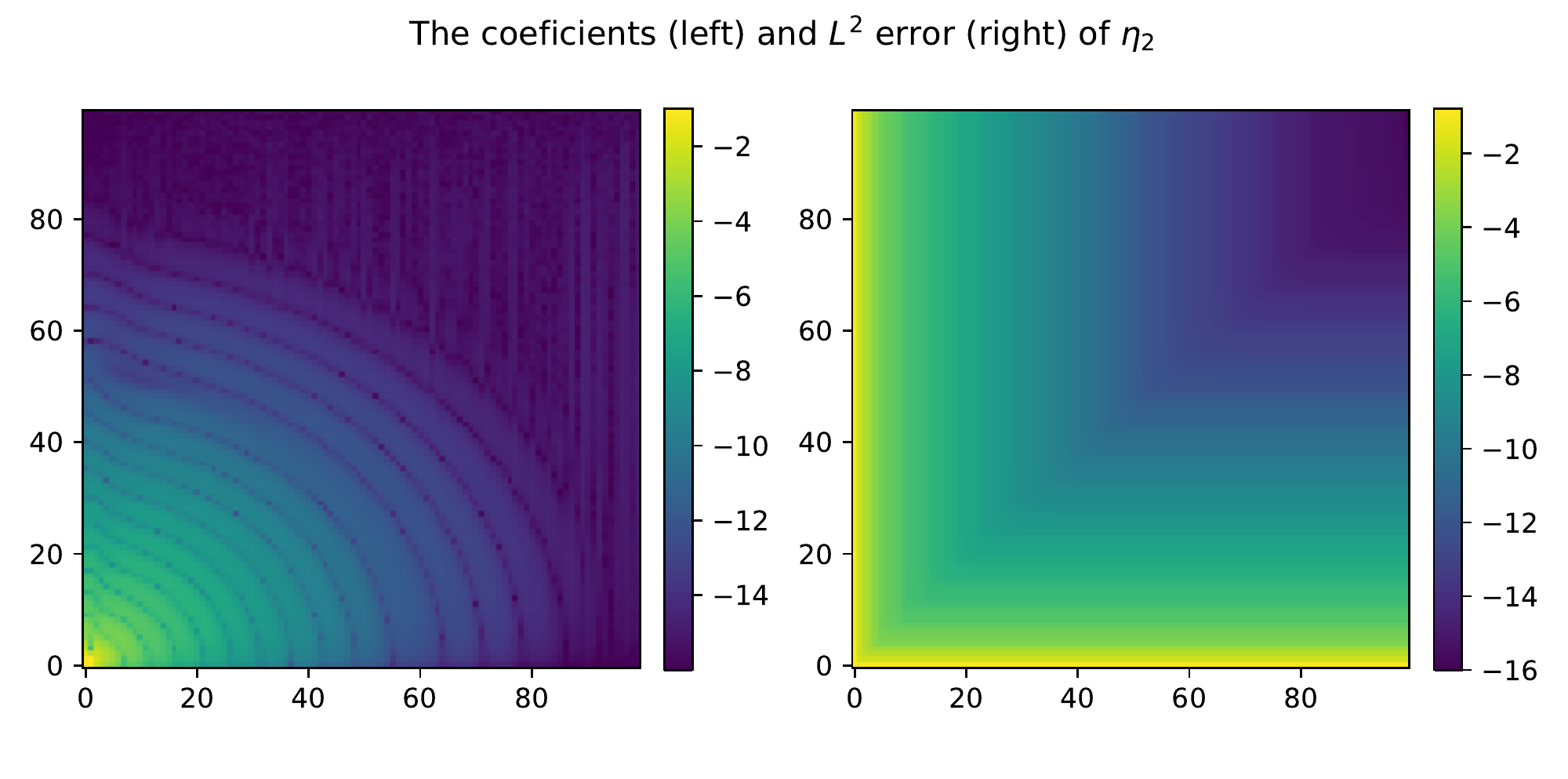}
	\includegraphics[width=\linewidth]{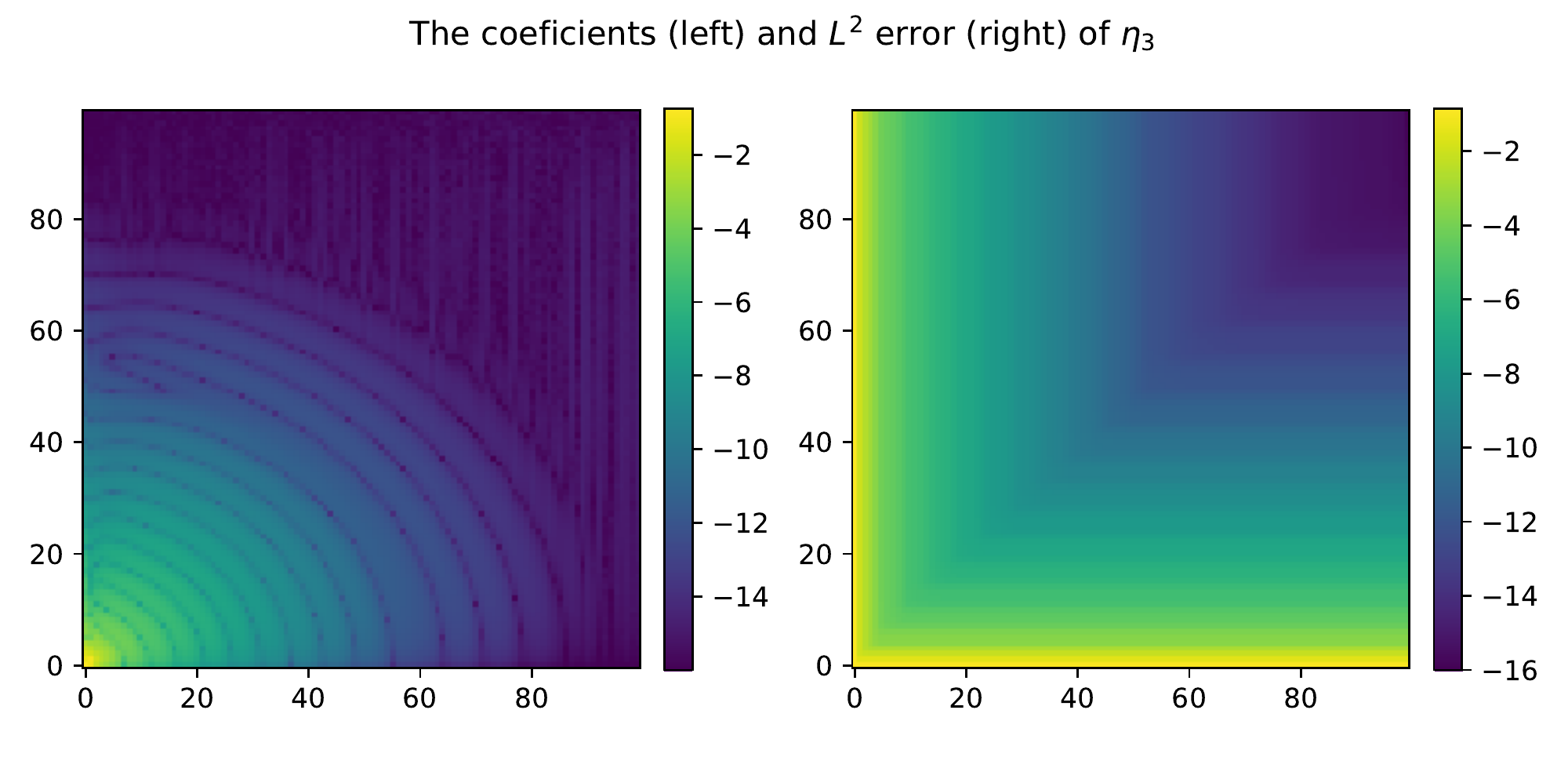}
	\caption{The coefficients (left) and the $L^2$ error of the Legendre
	approximation of $\eta_1$, $\eta_2$, $\eta_3$}
	\label{fig5.3}
\end{figure}

In Figure \ref{fig5.3}, we show the coefficients and $L^2$ error of the
Legendre approximations of $\eta_1$, $\eta_2$ and $\eta_3$, where
we take $n_1=100$, $n_2=100$ and $N_1=N_2=100$ in \eqref{eq5.11} and
\eqref{eq5.12}.
The results suggest that the Legendre expansion has spectral
accuracy and the error is about $10^{-15}$ when $n_1=90$, $n_2=75$.

After the coefficients $b_{st}^i$ of the Legendre approximation are
pre-calculated, we can efficiently calculate $\eta_i, i=1,2,3$ by
matrix-matrix multiplication for given values of $\mu_1, \mu_2$.  The other
elements in forth order moment tensor under diagonalized coordinates can be
obtained
by \eqref{eq:re_qandeta}. All the fourth order moments under original
coordinates
can be obtained by a coordinate transform.
The overall moment closure procedure is summarized in Algorithm
\ref{alg:repre_q_3d}.
\begin{algorithm}
	\setstretch{1.35}
	\caption{Calculate forth order moments for symmetric QEA on unit sphere}
	\label{alg:repre_q_3d}
	{\bf Input:}
	The values of the elements in second order moment $\hat{M}$  \\
	{\bf Output:}
	The values of the elements in fourth order moment $\hat{Q}$
	\begin{algorithmic}[1]
		\REQUIRE Pre-calculated coefficients $\{\, b_{st}^i, s=0,\ldots, n_1,
		t=0,\ldots, n_2, i=1,2,3\,\}$
		\STATE Use a linear algebra subroutine to calculate an orthogonal
		matrix $U$ which is formed by eigenvectors of $\hat{M}$ and make
        $M= U^T \hat{M} U$ into a diagonal matrix with $m_{11}\le m_{22} \le
        m_{33}$.
		\STATE Calculate the value of $\mu_1, \mu_2$ from $m_{11}, m_{22}$ by
		using \eqref{eq:m2mu}.
        \STATE Mapping the values of $\mu_1$ and $\mu_2$ to $x\times y \in
        [-1,1]\times[-1,1]$ with the transformation in Eq.(\ref{tr_mu_x})
		\STATE Calculate the values of $\eta_1$, $\eta_2$ and $\eta_3$ by
		Eq.(\ref{eq5.11})
		\STATE Calculate the nonzero elements in forth order moment under
		diagonalized coordinates by Eq.\eqref{eq:re_qandeta}.
        \STATE Use \eqref{eq:qhattrans} to calculate the elements in $\hat{Q}
        = (\hat{q}_{ijkl})$.
	\end{algorithmic}
\end{algorithm}

The overall computational time cost for evaluating $m$ moments in Algorithm
\ref{alg:repre_q_3d} is $O(n_1n_2m)$. Here, $n_1, n_2$ is about $70\sim
90$ to reach double precision. The storage cost is $O(n_1 n_2+ (n_1+n_2)m)$.
To reduce the computational time cost, we may use piecewise Legendre
approximation.
To show this, we divide the region of $(x,y) \in [-1,1]^2$ into
6 blocks.
These 6 blocks in $(x,y)$ are showed in  Figure \ref{fig:partirion-mu}(b)
marked by 6 different colors,  and the corresponding blocks in $(\mu_1,
\mu_2)$ are shown in Figure \ref{fig:partirion-mu}(a).
The coordinates of corner points in $(\mu_1, \mu_2)$ blocks are given
below:
\begin{align*}
&A_1(0,0), \quad A_2(\frac{1}{2},-\frac{1}{2}), \quad A_3(\frac{2}{3},0),
\quad B_1(\frac{1}{3},0), \quad B_2(\frac{5}{12},-\frac{1}{12}),\quad
B_3(\frac{7}{12},-\frac{1}{4}),  \\
&C_1(\frac{1}{9},0), \quad C_2(\frac{11}{72},-\frac{1}{24}), \quad
C_3(\frac{1}{4},-\frac{6}{35}), \quad C_4(\frac{19}{36},-\frac{5}{12}),\quad
D_1(\frac{1}{20},-\frac{1}{20}), \quad D_2(\frac{1}{6},-\frac{1}{6}).
\end{align*}
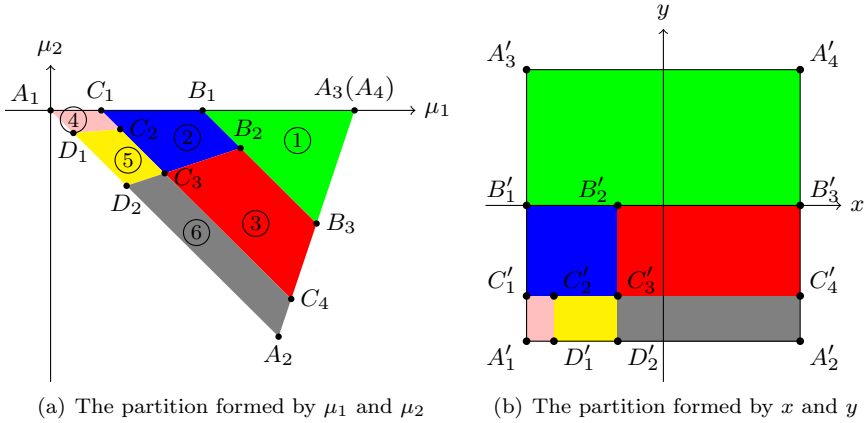
\begin{figure}[htbp]
	\centering{}
	\begin{subfigure}[The partition formed by $\mu_1$ and $\mu_2$]
		{\begin{tikzpicture}[scale=6]
		\fill[green] (2/3,0) -- (1/3,0) -- (7/12,-1/4);                              
		\fill[blue] (1/3,0) -- (5/12,-1/12) -- (1/4,-5/36) -- (1/9,0);               
		\fill[red] (7/12,-1/4) -- (5/12,-1/12) -- (1/4,-5/36) -- (19/36,-5/12);      
		\fill[pink] (1/9,0) -- (11/72,-1/24) -- (1/20,-1/20) -- (0,0);               
		\fill[yellow] (1/4,-5/36) -- (11/72,-1/24) -- (1/20,-1/20) -- (1/6,-1/6);    
		\fill[gray] (1/4,-5/36) -- (19/36,-5/12) -- (1/2,-1/2) -- (1/6,-1/6);        

		\node (a1) at (13/24,-1/15) {\circled{1}};
		\node (a2) at (0.3,-0.06) {\circled{2}};
		\node (a3) at (0.45,-0.25) {\circled{3}};
		\node (a4) at (0.05,-1/48) {\circled{4}};
		\node (a5) at (0.17,-0.11) {\circled{5}};
		\node (a6) at (0.32,-0.27) {\circled{6}};

		\fill (0,0) circle(0.2pt) node[above=6pt,left]{$A_1$};                       
		\fill (1/2,-1/2)circle(0.2pt) node[below]{$A_2$};
		\fill (2/3,0)circle(0.2pt) node[above]{$A_3$($A_4$)};
		\fill (1/3,0)circle(0.2pt) node[above]{$B_1$};
		\fill (5/12,-1/12)circle(0.2pt) node[right=3pt,above]{$B_2$};
		\fill (7/12,-1/4)circle(0.2pt) node[right]{$B_3$};
		\fill (1/9,0)circle(0.2pt) node[above]{$C_1$};
		\fill (11/72,-1/24)circle(0.2pt) node[right]{$C_2$};
		\fill (1/4,-5/36)circle(0.2pt) node[below=2pt,right]{$C_3$};
		\fill (19/36,-5/12)circle(0.2pt) node[right]{$C_4$};
		\fill (1/20,-1/20)circle(0.2pt) node[below]{$D_1$};
		\fill (1/6,-1/6)circle(0.2pt) node[below]{$D_2$};

		\draw[->] (-0.1,0) -- (0.8,0) node[right]{$\mu_1$};                          
		\draw[->] (0,-0.6) -- (0,0.1) node[above] {$\mu_2$};
		\end{tikzpicture}}
	\end{subfigure}
	\begin{subfigure}[The partition formed by $x$ and $y$]
		{\begin{tikzpicture}[scale=1.8]
			\fill[green] (-1,1) -- (-1,0) -- (1,0) -- (1,1);                              
			\fill[blue] (-1,0) -- (-1/3,0) -- (-1/3,-2/3) -- (-1,-2/3);                   
			\fill[red] (-1/3,0) -- (-1/3,-2/3) -- (1,-2/3) -- (1,0);                      
			\fill[pink] (-1,-2/3) -- (-4/5,-2/3) -- (-4/5,-1) -- (-1,-1);                 
			\fill[yellow] (-4/5,-1) -- (-1/3,-1) -- (-1/3,-2/3)-- (-4/5,-2/3);            
			\fill[gray] (-1/3,-2/3) -- (1,-2/3) -- (1,-1) -- (-1/3,-1);                   

			\fill (-1,-1) circle(0.8pt) node[below=6pt,left]{$A_1^{\prime}$};
			\fill (1,-1) circle(0.8pt) node[below=6pt,right]{$A_2^{\prime}$};
			\fill (-1,1) circle(0.8pt) node[above=6pt,left]{$A_3^{\prime}$};
			\fill (1,1) circle(0.8pt) node[above=6pt,right]{$A_4^{\prime}$};
			\fill (-1,0) circle(0.8pt) node[above=6pt,left]{$B_1^{\prime}$};
			\fill (-1/3,0) circle(0.8pt) node[above=6pt,left]{$B_2^{\prime}$};
			\fill (1,0) circle(0.8pt) node[above=6pt,right]{$B_3^{\prime}$};
			\fill (-1,-2/3) circle(0.8pt) node[above=6pt,left]{$C_1^{\prime}$};
			\fill (-4/5,-2/3) circle(0.8pt) node[above=6pt,right]{$C_2^{\prime}$};
			\fill (-1/3,-2/3) circle(0.8pt) node[above=6pt,right]{$C_3^{\prime}$};
			\fill (1,-2/3) circle(0.8pt) node[above=6pt,right]{$C_4^{\prime}$};
			\fill (-4/5,-1) circle(0.8pt) node[below=6pt,right]{$D_1^{\prime}$};
			\fill (-1/3,-1) circle(0.8pt) node[below=6pt,right]{$D_2^{\prime}$};

			\draw (-1,-1) -- (1,-1) -- (1,1) -- (-1,1) -- (-1,-1);
			\draw[->] (-1.3,0) -- (1.3,0) node[right] {$x$};
			\draw[->] (0,-1.3) -- (0,1.3) node[above] {$y$};
			\end{tikzpicture}}
	\end{subfigure}
	\caption{Partitions of parameter regions for piecewise Legendre
	approximations}
	\label{fig:partirion-mu}
\end{figure}

 We apply the Legendre approximation for each blocks in such a partition. The
 number of expansion terms $n_1$, $n_2$ and the corresponding $L^2$ error for
 $\eta_1$,
 $\eta_2$ and $\eta_3$ in each block are shown in Table \ref{table:1}. We see
 that $n_1$ and $n_2$ are both no
 more than 26 with the $L^2$ error gets below $10^{-14}$ in these blocks.
 Comparing to the global Legendre expansion in the whole region, the piecewise
 Legendre expansions greatly reduce the time cost of calculation.
\begin{table}
	\caption{The number of expansion terms $n_1$, $n_2$ and the $L^2$ error for
    $\eta_1$, $\eta_2$ and $\eta_3$ in each block with piecewise Legendre
    approximation.}
\label{table:1}
        \begin{tabular}{c|ccc|ccc|ccc}
		\hline
		\multirow{3}{*}{Block }
		&  \multicolumn{3}{|c|}{$\eta_1$}
		&  \multicolumn{3}{|c|}{$\eta_2$}
		&  \multicolumn{3}{|c}{$\eta_3$}
		\\ \cline{2-10}& $n_1$ & $n_2$ & $L^2$ error & $n_1$ & $n_2$ & $L^2$ error & $n_1$ & $n_2$ & $L^2$ error \\ \hline
		\circled{1}  & 17 & 15 & 7.26E-15 & 17 & 15 & 9.27E-15 & 17 & 16 & 9.07E-15\\ \hline
		\circled{2}  & 22 & 22 & 4.88E-15 & 22 & 21 & 6.48E-15 & 23 & 21 & 9.78E-15\\ \hline
		\circled{3}  & 16 & 20 & 6.11E-15 & 16 & 20 & 3.99E-15 & 16 & 21 & 8.18E-15\\ \hline
		\circled{4}  & 24 & 26 & 9.77E-15 & 23 & 26 & 8.63E-15 & 23 & 26 & 7.68E-15\\ \hline
		\circled{5}  & 23 & 24 & 9.71E-15 & 23 & 25 & 9.46E-15 & 22 & 26 & 8.04E-15\\ \hline
		\circled{6}  & 20 & 24 & 9.72E-15 & 21 & 24 & 6.23E-15 & 22 & 25 & 9.69E-15\\ \hline
	\end{tabular}
	\centering
\end{table}

\section{Summary}
We have shown some basic properties of quasi-equilibrium closure
approximation for antipodally symmetric problems and designed efficient
high order numerical implementations for such closure approximations on unit
circle and unit sphere by using global Legendre approximation and piecewise
Legendre approximation. The proposed implementation can reach to double
accuracy with mush smaller memory cost. The time efficiency is improved
by using piecewise polynomial approximations.
The proposed approach can be directly extended to handle other QEA closure
approximations, such as the Fisher-Bingham \cite{kent_fisherbingham_1982} and
von Mises-Fisher distributions\cite{sra_short_2012} for non-symmetric
problems.

Note that the tensor-product polynomial approximations are limited to
low-dimensional problems. For high-dimensional problems, spectral
sparse grid methods (see e.g. \cite{shen_efficient_2010}
\cite{shen_efficient_2012}) and deep neural
networks \cite{li_better_2020}\cite{yu_onsagernet_2020} are vital
approximation tools.
Implementation of high-dimensional QEA using these techniques will be the
topic of our future study.

\begin{acknowledgements}
The authors would like to thank Prof. Chuanju Xu, Li-Lian Wang
and Dr. Jie Xu  for helpful discussions.
This work is partially supported by
NNSFC Grant 11771439, 91852116 and China Science Challenge Project
no. TZ2018001.
\end{acknowledgements}

\noindent {\small {\bf Code and data availability}\
    All data and code generated or used during the
    study are available from the corresponding author by request.
}

\bibliographystyle{plain}  

\end{document}